\documentclass[11pt,leqno]{article}
\usepackage{amsthm,amsfonts,amssymb,epsfig,graphics,amsmath,
eufrak,oldgerm}
\usepackage[latin1]{inputenc}\relax
\usepackage[active]{srcltx}

\newcommand{\NN}{{\mathbb N}}

\newcommand{\RR}{{\mathbb R}}

\newcommand{\TT}{{\mathbb T}}

\newcommand{\ZZ}{{\mathbb Z}}

\def\cc{{\rm\bf c}}
\def\f{{\rm\bf f}}
\def\g{{\rm\bf g}}
\def\h{{\rm\bf h}}
\def\p{{\rm\bf p}}
\def\r{{\rm\bf r}}
\def\u{{\rm\bf u}}

\def\cE{\mathcal{E}}
\def\cF{\mathcal{F}}
\def\cG{\mathcal{G}}
\def\cH{\mathcal{H}}

\def\cJ{\mathcal{J}}

\def\cL{\mathcal{L}}

\def\cW{\mathcal{W}}

\def\cZ{\mathcal{Z}}

\def\mL{\mathfrak{L}}

\def\mP{\mathfrak{P}}

\def\mU{\mathfrak{U}}

\def\mX{\mathfrak{X}}

\def\ma{\mathfrak{a}}

\def\md{\mathfrak{d}}

\def\mg{\mathfrak{g}}
\def\mh{\mathfrak{h}}
\def\mi{\mathfrak{i}}

\def\mq{\mathfrak{q}}
\def\mr{\mathfrak{r}}

\def\mu{\mathfrak{u}}
\def\mv{\mathfrak{v}}

\def\eps {\varepsilon}

\def\part{\partial}

\def\u{{\bf u}}

\def\p{{\bf p}}
\def\r{{\bf r}}

\def\Dive{\md \mi \mv}
\def\Grade{\mg \mr \ma \md}

\newcommand\id{\mathrm{Id} }
\newcommand\Div{\mathrm{div }}

\newtheorem{theo}{Theorem} [section]
\newtheorem{lem}{Lemma} [section]

\newtheorem{prop}{Proposition} [section]

\numberwithin{equation}{section}

\begin{document} 

\centerline{\LARGE \bf Cascade of phases in turbulent flows}

\bigskip

\centerline{\sc Christophe CHEVERRY
\footnote{\footnotesize IRMAR, Universit\'e de Rennes I, Campus de 
Beaulieu, 35042 Rennes cedex, France, \text{ }\text{ }\text{ }\text{ }
christophe.cheverry@univ-rennes1.fr.}}

\bigskip
\bigskip

\noindent{\bf \small Abstract.} {\small This article is 
devoted to incompressible Euler equations (or to 
Navier-Stokes equations in the vanishing viscosity limit). 
It describes the propagation of {\it quasi-singularities}. 
The underlying phenomena are consistent with the notion 
of a {\it cascade of energy}.}

\bigskip

\noindent{\bf \small R\'esum\'e.} {\small Cet article 
\'etudie les \'equations d'Euler incompressible (ou de
Navier-Stokes en pr\'esence de viscosit\'e \'evanescente).
On y d\'ecrit la propagation de {\it quasi-singularit\'es}. 
Les ph\'enom\`enes sous-jacents confirment l'id\'ee selon 
laquelle il se produit une {\it cascade d'energie}.}

\bigskip

\setcounter{section}{0}

\section{Introduction.} $ \, $

\vskip -5mm

Consider incompressible fluid equations

\medskip

\noindent{$ (\cE) \qquad \part_t \u + (\u \cdot \nabla ) \, \u + 
\nabla \p = 0 \, , \qquad \Div \ \u = 0 \, , \qquad (t,x) \in 
[0,T] \times \RR^d \, , $}

\medskip

\noindent{where $ \u = {}^t ( \u^1 , \cdots , \u^d) $ is the fluid 
velocity  and $ \p $ is the pressure. The structure of {\it weak} 
solutions of $ (\cE) $ in $ d - $space dimensions with $ d \geq 2 $ 
is a problem of wide current interest \cite{Chem}-\cite{Li}. 
The questions are how to describe the phenomena with adequate 
models and how to visualize the results in spite of their
complexity. We will achieve a small step in these two 
directions. 

\medskip

According to the physical intuition, the appearance of 
singularities is linked with the {\it increase of the vorticity}.
Along this line, we have to mark the contributions \cite{BKM} 
and \cite{CF}. Interesting objects are solutions which do
not blow up in finite time but whose associated vorticities
increase arbitrarily fast. These are {\it quasi-singularities}.
Their study is of practical importance.

\medskip

Typical examples of quasi-singularities are oscillations.
This is a well-known fact going back to \cite{C}-\cite{MPP}.
The works \cite{C} and \cite{MPP} rely on phenomenological 
considerations and engineering experiments. Further 
developments are related to homogenization \cite{E}-\cite{E2}, 
compensated compactness \cite{D}-\cite{Ge} and non linear 
geometric optics \cite{Che}-\cite{CGM}-\cite{CGM1}. 

\smallskip

$\, $

\break

$ \, $

In \cite{D}, DiPerna and Majda show the persistence of 
oscillations in three dimensional Euler equations. To 
this end, they select a parameter $ \eps \in \, ]0,1] $ 
and look at
\begin{equation} \label{oscidipmaj}
\ \u^\eps_s (t,x) \, := \, {}^t \bigl( \g ( x_2 , \eps^{-1} 
\, x_2) , 0 , \h \bigl( x_1 - \g ( x_2 , \eps^{-1} \, x_2) 
\, t , x_2, \eps^{-1} \, x_2 \bigr) \bigr) \quad
\end{equation}
where $ \g (x_2,\theta) $ and $ \h (x_1,x_2,\theta) $ 
are smooth bounded functions with period $ 1 $ in 
$ \theta $. They remark that the functions $ \u^\eps_s $ 
are exact smooth solutions of $ (\cE) $ and they let $ \eps $
goes to zero. Yet, this construction is of a very special 
form. First, it comes from shear layers (steady 2-D solutions) as

\medskip

$ \tilde \u^\eps_s (t,x) \, = \, \tilde \u^\eps_s (0,x) \, 
= \, {}^t \bigl( \g ( x_2,\eps^{-1} \, x_2) , 0 \bigr) \, 
\in \, \RR^2 \, . $

\medskip

\noindent{Secondly, it involves a phase $ \varphi_0 (t,x) 
\equiv x_2 $ which does not depend on $ \eps $.} Of course,
this is a common fact \cite{CG}-\cite{G}-\cite{G2}-\cite{Se}
when dealing with large amplitude high frequency waves.
Nevertheless, this is far from giving a complete idea of what 
can happen.

\medskip

Our aim in this paper is to develop a theory which allows to 
remove the two restrictions mentioned above. Fix $ \flat = 
(l,N) \in \NN^2 $ where the integers $ l $ and $ N $ are such 
that $ 0 < l < N $. Introduce the {\it geometrical} phase

\medskip

$ \varphi^\eps_g (t,x) \, := \, \varphi_0 (t,x) \, + \, \sum_{k=1}^{
l-1} \, \eps^{\frac{k}{l}} \ \varphi_k (t,x) \, . $

\medskip

\noindent{In the section 2, we state the Theorem \ref{appBKW} which
provides with {\it approximate} solutions $ \u^\eps_\flat $ defined 
on the interval $ [0,T]  $ with $ T > 0 $ and having the form}
\begin{equation} \label{formegenerale}
\begin{array} {ll}
\u^\eps_\flat (t,x) \! \! \! & = \, {}^t ( \u^{\eps 1}_\flat, 
\cdots , \u^{\eps d}_\flat)(t,x) \\
\ & = \, \u_0 (t,x) \,  + \, \sum_{k=1}^N \, \eps^{
\frac{k}{l}} \ U_k \bigl( t,x, \eps^{-1} \, \varphi^\eps_g (t,x) 
\bigr) \qquad \qquad
\end{array}
\end{equation}
where the smooth profiles

\medskip

$ U_k (t,x,\theta) \, = \, {}^t (U^1_k, \cdots, U^d_k) 
(t,x,\theta) \in \RR^d \, , \qquad 1 \leq k \leq N \, , $

\medskip

\noindent{are periodic functions of $ \theta \in \RR /
\ZZ $.} We assume that 

\medskip

$ \exists \, (t,x,\theta) \in [0,T] \times \RR^d \times 
\TT \, ; \qquad \part_\theta U_1 (t,x,\theta) \not = 0 $.

\medskip

\noindent{We say that the family $ \{ \u^\eps_\flat \}_\eps $ 
is a {\it weak}, a {\it strong} or a {\it turbulent} 
oscillation according as we have respectively $ l=1 $, 
$ l = 2 $ or $ l \geq 3 $.}

\medskip

The order of magnitude of the energy of the oscillations 
is $ \eps^{\frac{1}{l}} $. Compute the vorticities associated with the 
functions $ \u^\eps_\flat $. These are the skew-symmetric 
matrices $ \Omega^\eps_\flat = (\Omega^{\eps i}_{\flat j} 
)_{1 \leq i,j \leq d} $ where
$$ \quad \begin{array} {l}
\Omega^{\eps i}_{\flat j} (t,x) \, := \, (\part_j 
\u^{\eps i}_\flat - \part_i \u^{\eps j}_\flat)(t,x) \\
\ = \, \sum_{k=1}^N \, \eps^{\frac{k}{l}-1} \ ( \part_j 
\varphi^\eps_g \ \part_\theta U^i_k - \part_i \varphi^\eps_g \ 
\part_\theta U^j_k) \bigl( t,x, \eps^{-1} \, \varphi^\eps_g (t,x) 
\bigr) \\
\quad \ + \, ( \part_j \u^i_0  - \part_i \u^j_0)(t,x) \, + \, 
\sum_{k=1}^N \, \eps^{\frac{k}{l}} \ ( \part_j U^i_k  - 
\part_i U^j_k ) \bigl( t,x, \eps^{-1} \, \varphi^\eps_g 
(t,x) \bigr) \, .  
\end{array} $$
The principal term in $ \Omega^\eps_\flat $ is of size 
$ \eps^{\frac{1}{l}-1}$. When $ l \geq 2 $, no uniform 
estimates are available on the family $ \{ \Omega^\eps_\flat 
\}_{\eps \in \, ]0,1] } $. In particular, if $ d = 3 $, there 
is no uniform control on the enstrophy

\medskip

$ \int_0^T \int_{\RR^3} \ \vert \omega^\eps_\flat (t,x) \vert^2 \ \, 
dt \, dx \, , \qquad \omega^\eps_\flat (t,x) \, := \, (\nabla \wedge 
\u^\eps_\flat)(t,x) \, \equiv \, \Omega^\eps_\flat (t,x) \, . $ 

\medskip

\noindent{We see here that strong and turbulent oscillations 
are examples of quasi-singularities.} Observe that the 
expansion (\ref{formegenerale}) involves a more complicated
structure than in (\ref{oscidipmaj}) though the corresponding
regime is less singular.

\medskip

The BKW analysis reveals that the phase shift $ \varphi_1 $
and the terms $ \varphi_k $ with $ 2 \leq k \leq l-1 $ play
different parts. The r\^ole of $ \varphi_1 $ is partly examined 
in the articles \cite{Che} and \cite{CGM} which deal with the 
case $ l = 2 $. When $ l \geq 3 $, the phenomenon to emphasize 
is the creation of the $ \varphi_k $ with $ 2 \leq k \leq l-1 $.
Indeed, suppose that

\medskip

$ \varphi_2 (0,\cdot) \equiv \, \cdots \, \equiv \varphi_{l-1} 
(0,\cdot) \equiv 0 \, , \qquad l \geq 3 \, . $

\medskip

\noindent{Then, generically, we find}

\medskip

$ \exists \, t \in \, ]0,T] \, ; \qquad \varphi_2 (t,\cdot) 
\not \equiv 0 \, , \qquad \cdots \, , \qquad \varphi_{l-1} 
(t,\cdot) \not \equiv 0 \, . $

\medskip

\noindent{Now starting with {\it large} amplitude waves (this 
corresponds to the limit case $ l = + \infty $) that is}

\medskip

$ \u^\eps_\infty (0,x) \, = \, \sum_{k=0}^\infty \, 
\eps^k \ U_k \bigl( 0,x, \eps^{-1} \, \varphi_0 (0,x) 
\bigr) \, , \qquad  \part_\theta U_0 \not \equiv 0 \, , $

\medskip

\noindent{the description of $ \u^\eps_\infty (t,\cdot) $ on the 
interval $ [0,T] $ with $ T > 0 $ needs the introduction of an {\it
infinite cascade} of phases $ \varphi_k $.} The scenario is 
the following. Oscillations of the velocity develop spontaneously 
in all the intermediate frequencies $ \eps^{\frac{k}{l}-1} $ 
and in all the directions $ \nabla \varphi_k (t,x) $. This 
expresses {\it turbulent} features in the flow.

\medskip

The family $ \{ \u^\eps_\flat \}_{\eps \in \, ]0,1]} $
is $ \eps - $stratified \cite{G2} with respect to the phase
$ \varphi^\eps_g $ with in general $ \varphi^\eps_g \not 
\equiv \varphi_0 $. The presence in  $ \varphi^\eps_g $ of the
non trivial functions $ \varphi_k $ is necessary and sufficient 
to encompass the {\it geometrical} features of the propagation. 
It has various consequences which are detailed in the section 
3. It brings informations about microstructures, compensated 
compactness and non linear geometric optics. It also confirms 
observations made in the statistical approach of turbulences 
\cite{FMRT}-\cite{L}.

\medskip

The chapter 4 is devoted to the demonstration of Theorem
\ref{appBKW}. Because of {\it closure problems}, the use of 
the geometrical phase $ \varphi^\eps_g $ does not suffice to 
perform the BKW analysis. Among other things,  {\it adjusting 
phases} $ \varphi_k $ with $ l \leq k \leq N $ must be 
incorporated in order to put the system of formal equations 
in a triangular form.

\medskip

The expressions $ \u^\eps_\flat $ are not exact solutions of
Euler equations, yielding small error terms $ \f^\eps_\flat $ 
as source terms. The matter is to know if there exists exact 
solutions which coincide with $ \u^\eps_\flat (0,\cdot) $ at 
time $ t = 0 $, which are defined on $ [0,T] $ with $ T > 0 $, 
and which are close to the approximate divergence free solutions 
$ \u^\eps_\flat $. This is a problem of {\it stability}.

\medskip

\noindent{The construction of exact solutions requires a good 
understanding of the different mechanisms of amplifications 
which occur.} In the subsection 5.1, we make a distinction 
between {\it  obvious} and {\it hidden} instabilities.

\medskip

\noindent{The obvious instabilities can be detected by 
looking at the BKW analysis presented before.} They
imply the non linear instability of Euler equations
(Proposition \ref{euinstap}). They need to be absorbed 
a dependent change of variables which induces a defect 
of hyperbolicity. The hidden instabilities can be revealed 
by soliciting this lack of hyperbolicity. They require to 
be controled the addition of dissipation terms.

\medskip

\noindent{In the subsection 5.2, we look at incompressible fluids 
with anisotropic viscosity.} This is the framework of \cite{CDGG}
though we adopt a different point of view. We consider strong 
oscillations. We show (Theorem \ref{ciprin}) that {\it exact}
solutions corresponding to $ \u^\eps_{(2,N)} $ exist on some 
interval $ [0,T] $ with $ T > 0 $ independent on $ \eps \in 
\, ]0,1] $. 

\bigskip
\smallskip

{\small \parskip=-3pt
\tableofcontents
}

\break

\section{Euler equations in the variables $ (t,x) $.} $\,$

\vskip -5mm

The description of incompressible flows in turbulent regime 
is a delicate question. No systematic analysis is yet available. 
However, special appro- ximate solutions with rapidly varying 
structure in space and time can be exhibited. Their construction
is summarized in this chapter 2.

\subsection{Notations.}

\smallskip
  
\noindent{$ \bullet $ {\em Variables.}} Let $ T \in \RR^+_* $. The time 
variable is $ t \in [0,T] $. Let $ d \in \NN \setminus \{ 0 ,
1 \} $. The space variables are $ (x, \theta) \in \RR^d \times
\TT $ where $ \TT := \RR / \ZZ $. Mark the ball

\medskip

$ B(0,R] \, := \, \bigl \{ \, x \in \RR^d \, ; \ \vert x \vert^2
:= \sum_{i=1}^d \, x_i^{\, 2} \leq R \, \bigr \} \, , \qquad R \in 
\RR^+ \, . $ 

\medskip

\noindent{The state variables are the velocity field $ u = {}^t 
(u^1,\cdots,u^d) \in \RR^d $ and the pressure $ p \in \RR $. 
Given $ (u , \tilde u ) \in (\RR^d)^2 $, define

\medskip

$ u \cdot \tilde u := \sum_{i=1}^d \, u^i \, \tilde u^i \, , 
\qquad \vert u \vert^2 := u \cdot u \, , \qquad u \otimes \tilde 
u := (u^j \, \tilde u^i )_{1 \leq i,j \leq d} \, . $

\medskip

\noindent{The symbol $ S^d_+ $ is for the set of positive definite 
quadratic form on $ \RR^d $.} An element $ \mq \in S^d_+ $ can be 
represented by some $ d \times d $ matrix $ ( \mq_{ij} )_{1 \leq 
i,j \leq d} $.

\bigskip

\noindent{$ \bullet $ {\em Functional spaces.}} Distinguish the 
expressions $ \u (t,x) $ which do not depend on the variable 
$ \theta $ from the expressions $ u(t,x,\theta) $ which depend on 
$ \theta $. The boldfaced type $ \u $ is used in the first
case whereas the letter $ u $ is employed in the second 
situation. 

\medskip

\noindent{Note $ C^\infty_b ( [0,T] \times \RR^d ) $ the space 
of functions in $ [0,T] \times \RR^d $ with bounded continuous 
derivatives of any order.} Let $ m \in \NN $. The Sobolev space 
$ H^m $ is the set of functions 

\medskip

$ u(x,\theta) \, = \, \sum_{k \in \ZZ} \, \u_k (x) \ e^{
i \, k \, \theta} $ 

\medskip

\noindent{such that}

\medskip

$ \parallel u \parallel_{H^m}^2 \, := \, \sum_{k \in \ZZ} \, (1 +
\vert k \vert^2)^m \ \int_{\RR^d} \, (1 + \vert \xi \vert^2)^m \
\vert \hat \u_k (\xi) \vert^2 \ d \xi \, < \, \infty $

\medskip

\noindent{where} 

\medskip

$ \cF ( \u )(\xi) \, = \, \hat \u (\xi) \, := \, ( 2 \, \pi)^{- 
\frac{d}{2}} \ \int_{\RR^d} \, e^{- i \, x \cdot \xi} \ \u (x) \ 
dx \, , \qquad \xi \in \RR^d \, . $

\medskip

\noindent{With these conventions, the condition $ \u \in H^m $ 
means simply that}

\medskip

$ \parallel \u \parallel_{H^m}^2 \, := \, \int_{\RR^d} \, (1 + 
\vert \xi \vert^2)^m \ \vert \hat \u (\xi) \vert^2 \ d \xi \, 
< \, \infty \, . $

\medskip

\noindent{Define}

\medskip

$ H^m_T \, := \, \bigl \{ \, u \ ; \ \part^j_t u \in L^2 ( [0,T];
H^{m-j}) \, , \ \forall \, j \in \{0, \cdots,m \} \, \bigr \} \, , $

\medskip

$ \cW^m_T \, := \, \bigl \{ \, u \ ; \ u \in C^j ( [0,T]; H^{m-j}) \, , 
\ \forall \, j \in \{0, \cdots,m \} \, \bigr \} \, , $

\medskip

\noindent{with the corresponding norms}

\medskip

$ \parallel u \parallel_{H^m_T}^2 \, := \, \sum_{j=0}^m \, \int_0^T \,
\parallel \part^j_t u (t,\cdot) \parallel_{H^m}^2 \ dt \, , $

\medskip

$ \parallel u \parallel_{\cW^m_T} \, := \, \sup_{t \in [0,T]} \ \ 
\sum_{j=0}^m \, \parallel \part^j_t u (t,\cdot) \parallel_{H^m} \, . $

\medskip

\noindent{Consider also }
$$ \ \left. \begin{array} {lll}
H^m_\infty \, := \, \bigcap_{T \in \RR^+} \, H^m_T \, , & \quad 
H^\infty_T \, := \, \bigcap_{m \in \NN} \, H^m_T \, , & \quad 
H^\infty_\infty \, := \, \bigcap_{T \in \RR^+} \, H^\infty_T \, , \\ 
\cW^m_\infty \, := \, \bigcap_{T \in \RR^+} \, \cW^m_T \, , & \quad  
\cW^\infty_T \, := \, \bigcap_{m \in \NN} \, \cW^m_T \, , & \quad 
\cW^\infty_\infty \, := \, \bigcap_{T \in \RR^+} \, \cW^\infty_T \, . 
\end{array} \right. $$
When $ m = 0 $, replace $ H^0 $ with $ L^2 $. Any function $ u 
\in L^2 $ can be decomposed according to

\medskip

$ u (t,x,\theta) \, = \, \langle u \rangle (t,x) + u^* (t,x,\theta) \, 
= \, \bar u (t,x) + u^* (t,x,\theta) $

\medskip

\noindent{where}

\medskip

$ \langle u \rangle (t,x) \, \equiv \, \bar u (t,x) \, := \, \int_\TT \, 
u(t,x,\theta) \ d \theta \, . $ 

\medskip

\noindent{Let $ \Gamma $ be the symbol of any of the spaces $  H^m $, 
$ H^m_T $, $ \cW^m_T $, $ \cdots $ defined before.} In order to
specify the functions with mean value zero, introduce

\medskip

$ \Gamma^* \, := \, \lbrace \, u \in \Gamma \, ; \ \bar u \equiv 0 \, 
\rbrace \, . $

\medskip

\noindent{Mark also}

\medskip

$ \text{supp}_x \, u^* \, := \, \text{closure of} \ \bigl \{ \, x 
\in \RR^d \, ; \ \parallel u^* (x, \cdot) \parallel_{L^2(\TT)} 
\not = \, 0 \, \bigr \} \, . $

\bigskip

\noindent{$ \bullet $ {\em Differential operators.}} Note

\medskip

$ \part_t \equiv \part_0 := \part / \part \, t \, , \qquad \ \
\part_\theta \equiv \part_{d+1} := \part / \part \, \theta \, , $

\smallskip

$ \part_j := \part / \part \, x_j \, , \qquad \qquad \ j 
\in \{1 , \cdots , d \} \, , $

\smallskip

$ \nabla := (\part_1, \cdots ,\part_d) \, , \qquad \Delta := 
\Delta_x + \part^2_\theta = \part^{\, 2}_1 + \cdots + \part^{\, 
2}_{d} + \part^2_\theta \, . $

\medskip

\noindent{Let $ u \in \cW^\infty_T $.} Define

\medskip

$ u \cdot \nabla \, := \, u^1 \ \part_1 + \cdots + u^d \ \part_d \, , $

\smallskip

$ \Div \ u \, := \, \part_1 u^1 + \cdots + \part_d u^d \, , $

\smallskip

$ \Div \ (u \otimes \tilde u) \, := \, \sum_{j=1}^d \, {}^t \bigl( 
\part_j ( u^j \, \tilde u^1) \, , \cdots , \part_j ( u^j \, \tilde 
u^d) \bigr) \in \RR^d \, . $

\medskip

\noindent{Employ the bracket $ < \cdot , \cdot >_H $ for the scalar 
product in the Hilbert space $ H $. Note $ \cL (E;F) $ the space of 
linear continuous applications $ T : E \longrightarrow F $ where $ E $ 
and $ F $ are Banach spaces. The symbol $ \cL (E) $ is simply for 
$\cL (E;E) $. Introduce the commutator

\medskip

$ [A;B] \, := \, A \circ B - B \circ A \, , \qquad (A,B) \in \cL(E)^2 \, . $

\medskip

\noindent{Let $ r \in \ZZ $.} The operator $ T $ is in $ \mL^r $ if

\medskip

$ \parallel T \parallel_{\cL (H^{m+r}_T ; H^m_T)} \, < \, \infty \, , 
\qquad \forall \, m \in \NN \, . $

\medskip

\noindent{Let $ \eps_0 > 0 $.} The family of operators $ \{ T^\eps 
\}_\eps \in \cL (H^\infty_T)^{]0,\eps_0]} $ is in $ \mU \mL^r $ if 

\medskip

$ \sup_{\eps \in \, ]0,\eps_0]} \ \ \parallel T^\eps \parallel_{
\cL (H^{m+r}_T ; H^m_T)} \, < \, \infty \, , \qquad \forall \, m 
\in \NN \, . $

\medskip

\noindent{Consider a family $ \{ f^\eps \}_\eps \in (\cW^\infty_T)^{
]0,\eps_0]} $. We say that $ \{ f^\eps \}_\eps = \bigcirc (\eps^r ) $ 
if

\medskip

$ \sup_{\eps \in \, ]0,\eps_0]} \ \ \eps^{-r} \ \parallel f^\eps 
\parallel_{\cW^m_T} \, < \, \infty \, , \qquad \forall \, m \in 
\NN \, . $

\medskip

\noindent{Given a family $ \{ \f^\eps \}_\eps \in (\cW^\infty_T)^{
]0,\eps_0]} $, we say that $ \{ \f^\eps \}_\eps = \bigcirc (\eps^r ) $ 
if

\medskip

$ \sup_{\eps \in \, ]0,\eps_0]} \ \ \eps^{-r+m} \ \parallel \f^\eps 
\parallel_{\cW^m_T} \, < \, \infty \, , \qquad \forall \, m \in 
\NN \, . $

\medskip

\noindent{Observe that the two preceding definitions have very 
different significations according as we use the letter $ f $ 
or the boldfaced type $ \f $.} In particular, the second inequalities
correspond to $ \eps - \, $stratified estimates. The families 
$ \{ f^\eps \}_\eps $ or $ \{ \f^\eps \}_\eps $ are $ \bigcirc 
(\eps^\infty ) $ if they are $ \bigcirc (\eps^r ) $ for all 
$ r \in \RR $.

\subsection{Divergence free approximate solutions in $ (t,x) $.} 

\noindent{$ \bullet $ {\bf A first result.}} Select smooth functions

\medskip

$ \u_{00} \in H^\infty \, , \qquad \varphi_{00} \in C^1 (\RR^d ) \, ,
\qquad \nabla \varphi_{00} \in C^\infty_b (\RR^d ) \, . $

\medskip

\noindent{Suppose that} 

\medskip

$ \exists \, c > 0 \, ; \qquad \vert \nabla \varphi_{00} (x) 
\vert \, \geq \, 2 \ c \, , \qquad \forall \, x \in \RR^d \, . $

\medskip

\noindent{For $ T > 0 $ small enough, the equation $ (\cE) $
associated with}

\medskip

$ \u_0 (0,x) = \u_{00} (x) \, , \qquad \forall \, x \in \RR^d $}

\medskip

\noindent{has a smooth solution $ \u_0(t,x) \in \cW^\infty_T $.}
Solve the eiconal equation

\medskip

\noindent{$ (ei) \qquad \part_t \varphi_0 + (\u_0 \cdot \nabla) \, 
\varphi_0 = 0 \, , \qquad (t,x) \in [0,T] \times \RR^d $}

\medskip

\noindent{with the initial data}

\medskip

$ \varphi_0 (0,x) = \varphi_{00} (x) \, , \qquad \forall \, x \in 
\RR^d \, . $}

\medskip

\noindent{If necessary, restrict the time $ T $ in order to have}
\begin{equation} \label{nonstaou} 
\vert \nabla \varphi_0 (t,x) \vert \, \geq \, c \, , \qquad \forall \, 
(t,x) \in [0,T] \times \RR^d \, . \qquad \qquad \qquad
\end{equation}
Call $ \Pi_0 (t,x) $ the orthogonal projector from $ \RR^d $ 
onto the hyperplane

\medskip

$ \nabla \varphi_0 (t,x)^\perp \, := \, \bigl \lbrace \, u 
\in \RR^d \, ; \ u \cdot \nabla \varphi_0 (t,x) = 0 \, \bigr 
\rbrace \, .$}

\smallskip  

\begin{theo} \label{appBKW} Select any $ \flat = (l,N) \in \NN^2_* $ 
such that $ 0 < l \, (3 + \frac{d}{2}) \ll N $. Consider the following 
initial data

\medskip

$ U_{k0}^* (x,\theta) = \Pi_0 (0,x) \, U_{k0}^* (x,\theta) \in 
H^\infty \, , \qquad  1 \leq k \leq N \, , $

\smallskip

$ \bar U_{k0} (x) \in H^\infty \, , \qquad 1 \leq k \leq 
N \, , $

\smallskip

$ \varphi_{k0} (x) \in H^\infty \, , \qquad 1 \leq k \leq l-1 \, . $

\medskip

\noindent{First, there are finite sequences $ \{ U_k \}_{1 \leq k \leq N} $ 
and $ \{ P_k \}_{1 \leq k \leq N} $ with}

\medskip

$ U_k (t,x,\theta) \in \cW^\infty_T \, , \qquad P_k (t,x,\theta) \in 
\cW^\infty_T \, , \qquad 1 \leq k \leq N \, , $

\medskip

\noindent{and a finite sequence $ \{ \varphi_k \}_{1 \leq k \leq l-1} $ 
with}

\medskip

$ \varphi_k (t,x) \in \cW^\infty_T \, , \qquad 1 \leq k \leq l-1 \, , $

\medskip

\noindent{which are such that}

\medskip

$ \Pi_0 (0,x) \, U_k^* (0,x,\theta) =  \Pi_0 (0,x) \, U_{k0}^*
(x,\theta) \, , \qquad 1 \leq k \leq N \, , $

\smallskip

$ \bar U_k (0,x) = \bar U_{k0} (x) \, , \qquad \! 1 \leq k \leq 
N \, , $

\smallskip

$ \varphi_k (0,x) = \varphi_{k0} (x) \, , \qquad \! 1 \leq k 
\leq l-1 \, . $ 

\medskip

\noindent{Secondly, there is $ \eps_0 \in \, ]0,1 ] $ and  
correctors}

\medskip

$ \cc \u^\eps_\flat (t,x) \in \cW^\infty_T \, , \qquad \cc \p^\eps_\flat 
(t,x) \in \cW^\infty_T \, , \qquad \eps \in \,  ]0,\eps_0] \, , $ 

\medskip

\noindent{which give rise to families satisfying}

\medskip

$ \{ \cc \u^\eps_\flat \}_\eps \, = \, \bigcirc(\eps^{\frac{N}{l}-2}) \, , 
\qquad \{ \cc \p^\eps_\flat \}_\eps \, = \, \bigcirc(\eps^{\frac{N}{l}}) 
\, . $ 

\medskip

\noindent{Then, all these expressions are adjusted so that the functions 
$ \u^\eps_\flat $ and $ \p^\eps_\flat $ defined according to}
\vskip -4mm
\begin{equation} \label{BKWdel}
\, \left. \begin{array}{l}
\u^\eps_\flat (t,x) := \u_0 (t,x) + \sum_{k=1}^N \, 
\eps^{\frac{k}{l}} \ U_k \bigl( t,x, \eps^{-1} \, \varphi^\eps_g 
(t,x) \bigr) + \cc \u^\eps_\flat (t,x) \ \, \\
\p^\eps_\flat (t,x) := \p_0 (t,x) + \sum_{k=1}^N \, 
\eps^{\frac{k}{l}} \ P_k \bigl( t,x, \eps^{-1} \, \varphi^\eps_g 
(t,x) \bigr) + \cc \p^\eps_\flat (t,x) \ \ \,
\end{array} \right.
\end{equation}
where $ \varphi^\eps_g (t,x) $ is the geometrical phase 
\vskip -4mm  
\begin{equation} \label{phasegeo} 
\varphi^\eps_g (t,x) \, := \, \varphi_0 (t,x) \, + \, \hbox{$ \sum_{k=1}^{l-1}$}
\ \eps^{\frac{k}{l}} \ \varphi_k (t,x) \qquad \qquad \qquad \ \
\end{equation}
are approximate solutions of $ (\cE) $ on the interval $ [0,T] $.
More precisely 

\medskip

$ \part_t \u^\eps_\flat + ( \u^\eps_\flat \cdot \nabla) \u^\eps_\flat 
+ \nabla \p^\eps_\flat = \f^\eps_\flat \, , \qquad \Div \ \u^\eps_\flat 
= 0 \, , \qquad \f^\eps_\flat = \bigcirc(\eps^{\frac{N}{l}-3-\frac{d}
{2}}) \, . $

\end{theo}

\smallskip

\noindent{$ \bullet $ {\bf Some comments.}} 

\medskip

\noindent{{\it Remark 2.2.1:}} In what follows, we suppose that 
$ U^*_1 $ is non trivial.} In other words, we start with some 
initial data satisfying
\begin{equation} \label{nonnontri}
\ \exists \, (x,\theta) \in \RR^d \times \TT \, ; \qquad U^*_1
(0,x,\theta) = U^*_{10} (x,\theta) \not = 0 \, . \qquad \qquad 
\qquad \quad
\end{equation}

\vskip -9mm

\hfill $ \triangle $

\bigskip

\noindent{{\it Remark 2.2.2:}} Fix any $ l \in \NN_* $. The Borel's
summation process allows to take $ N = + \infty $ in the Theorem
\ref{appBKW}. It yields BKW solutions $ (\u^\eps_\flat , \p^\eps_\flat) $ 
which solve $ (\cE) $ with infinite accuracy

\medskip

$ \part_t \u^\eps_\flat + ( \u^\eps_\flat \cdot \nabla) \u^\eps_\flat + \nabla 
\p^\eps_\flat \, = \, \bigcirc(\eps^\infty) \, , \qquad \Div \ \u^\eps_\flat 
\, = \, 0 \, . $ \hfill $ \triangle $

\bigskip

\noindent{{\it Remark 2.2.3:}} Suppose that the function $ \u_0 \in 
\cW^\infty_\infty $ is a global solution of Euler equations. Suppose 
also that the phase $ \varphi_0 \in \cW^\infty_\infty $ is subjected 
to (\ref{nonstaou}) on the strip $ [0,\infty [ \times \RR^d $ and that 
it is a global solution of the eiconal equation $ (ei) $. Then 
the Theorem \ref{appBKW} can be applied with any $ T \in \RR^+_* $. It 
means that no blow up occurs at the level of the equations yielding the 
profiles $ U_k $, $ P_k $ and the phases $ \varphi_k $. Yet, non linear
effects are present. \hfill $ \triangle $

\bigskip

\noindent{{\it Remark 2.2.4:}} The characteristic curves
of the field $ \part_t + \u_0 \cdot \nabla_x $ are obtained
by solving the differential equation 

\medskip

$ \part_t \, \Gamma (t,x) \, = \, \u_0 \bigl( t, \Gamma (t,x)
\bigr) \, , \qquad \Gamma (0,x) = x \, . $

\medskip

\noindent{Suppose that the oscillations of the profiles 
$ U^*_{k0} $ are concentrated in some domain $ D \subset 
\RR^d $.} In other words 

\medskip

$ \text{supp}_x \, U_{k0}^* \, \subset \, D \, , \qquad 
\forall \, k \in \{ 1, \cdots , N \} \, . $

\medskip

\noindent{The BKW analysis reveals that for all $ t \in 
[0,T] $ we have} 

\medskip

$ \text{supp}_x \, U_k^* (t,\cdot) \, \subset \, \bigl 
\lbrace \, \Gamma (t,x) \, ; \ x \in D \, \bigr \rbrace 
\, , \qquad \forall \, k \in \{ 1, \cdots , N \} \, . $

\medskip

\noindent{The phenomena under study have a finite speed of 
propagation.} \hfill $ \triangle $

\bigskip

\noindent{{\it Remark 2.2.5:}} The influence of dissipation 
terms will be taken into account in the subsection 4.1. The 
viscosity we will incorporate is anisotropic. It is small enough
in the direction $ \nabla \varphi^\eps_\flat $ in order to be 
compatible with the propagation of oscillations.  \hfill $ \triangle $ 

\subsection{End of the proof of Theorem \ref{appBKW}.} 

The Theorem \ref{appBKW} is a consequence of the Proposition 
\ref{appBKWinter} which will be stated and demonstrated in the 
subsection 4.2. Below, we just explain how to deduce the Theorem 
\ref{appBKW} from the Proposition \ref{appBKWinter} applied with 
$ \nu = 0 $. 

\medskip

\noindent{$ \bullet $ {\bf Dictionary between the profiles}.} Select
arbitrary initial data for

\medskip

$ \Pi_0 (0,x) \, \tilde U_k^* (0,x,\theta) \in H^\infty \, , \qquad 
\langle \tilde U_k \rangle (0,x) \in H^\infty \, , \qquad 1 \leq k
\leq N \, , $

\medskip

\noindent{and arbitrary initial data for}

\medskip

$ \varphi_k (0,x) \in H^\infty \, , \qquad  1 \leq k \leq l-1 \, . $

\medskip

\noindent{On the contrary, impose}
\begin{equation} \label{inipreli}
\varphi_k (0,\cdot) \equiv 0 \, , \qquad \forall \, k \in \{l,
\cdots, N \} \, . \qquad \qquad \qquad \qquad \qquad
\end{equation}
The Proposition \ref{appBKWinter} provides with finite 
sequences

\medskip 

$ \{ \tilde U_k \}_{1 \leq k \leq N} \, , \qquad \{ \tilde P_k \}_{
1 \leq k \leq N} \, , \qquad \{ \varphi_k \}_{1 \leq k \leq N} \, , $

\medskip

\noindent{and source terms}

\medskip

$ \tilde f^\eps_\flat (t,x,\theta) \in \cW^\infty_T \, , \qquad 
\tilde g^\eps_\flat (t,x,\theta) \in \cW^\infty_T \, . $

\medskip

\noindent{such that the associated oscillations}
$$ \left. \begin{array}{l}
\tilde \u^\eps_\flat (t,x) \, := \, \u_0 (t,x) \, + \, 
\sum_{k=1}^{N} \, \eps^{\frac{k}{l}} \ \tilde U_k \bigl( 
t,x, \eps^{-1} \, \varphi^\eps_\flat (t,x) \bigr) \, , \\
\tilde \p^\eps_\flat (t,x) \, := \, \p_0 (t,x) \, + \, 
\sum_{k=1}^{N} \, \eps^{\frac{k}{l}} \ \tilde P_k \bigl(t,
x, \eps^{-1} \, \varphi^\eps_\flat (t,x) \bigr) \, , \qquad 
\qquad \qquad \ \\
\tilde \f^\eps_\flat (t,x) \, := \, \eps^{-1} \ \tilde 
f^\eps_\flat \bigl(t,x, \eps^{-1} \, \varphi^\eps_\flat 
(t,x) \bigr) \, , \\
\tilde \g^\eps_\flat (t,x) \, := \, \eps^{-1} \ \tilde 
g^\eps_\flat \bigl(t,x, \eps^{-1} \, \varphi^\eps_\flat 
(t,x) \bigr) \, ,
\end{array} \right. $$
are subjected to

\medskip

$ \part_t \tilde \u^\eps_\flat + ( \tilde \u^\eps_\flat \cdot \nabla) 
\tilde \u^\eps_\flat + \nabla \tilde \p^\eps_\flat = \tilde \f^\eps_\flat
= \bigcirc ( \eps^{\frac{N+1}{l}-1}) \, , \quad \ \Div \ \tilde \u^\eps_\flat 
= \tilde \g^\eps_\flat = \bigcirc ( \eps^{\frac{N+1}{l}-1}) \, . $

\medskip

\noindent{The oscillations $ \tilde \u^\eps_\flat $ and $ \tilde \p^\eps_\flat $
involve the {\it complete} phase $ \varphi^\eps_\flat (t,x) $ which is
the sum of the geometrical phase $ \varphi^\eps_g (t,x) $ plus some 
{\it adjusting} phase $ \eps \, \varphi^\eps_a (t,x) $. More precisely

\medskip

$ \varphi^\eps_\flat (t,x) := \varphi^\eps_g (t,x) + \eps \ 
\varphi^\eps_a (t,x) \, , \qquad \varphi^\eps_a (t,x) := 
\sum_{k=l}^N \, \eps^{\frac{k}{l}-1} \ \varphi_k (t,x) \, . $

\medskip

\noindent{The functions $ \tilde \u^\eps_\flat $ and $ \tilde 
\p^\eps_\flat $ can also be written in terms of the phase 
$ \varphi^\eps_g $.} Indeed, there is a unique decomposition

\medskip

$ \tilde \u^\eps_\flat = \u^\eps_\flat + \r \u^\eps_\flat =
\u^\eps_\flat + \bigcirc (\eps^{\frac{N+1}{l}}) \, , \qquad 
\tilde \p^\eps_\flat = \p^\eps_\flat + \r \p^\eps_\flat =
\p^\eps_\flat + \bigcirc (\eps^{\frac{N+1}{l}}) \, , $

\medskip

\noindent{involving the representations}
\begin{equation} \label{BKWdelbis}
\ \u^\eps_\flat (t,x) = u^\eps_\flat \bigl( t,x, \eps^{-1} \, 
\varphi^\eps_g (t,x) \bigr) \, , \qquad \p^\eps_\flat (t,x) = 
p^\eps_\flat \bigl( t,x, \eps^{-1} \, \varphi^\eps_g (t,x) 
\bigr) 
\end{equation}
where the profiles $ u^\eps_\flat ( t,x, \theta ) $ and 
$ p^\eps_\flat ( t,x, \theta ) $ have the form
$$ \left. \begin{array}{l}
u^\eps_\flat ( t,x, \theta ) \, = \, \u_0 (t,x) + \sum_{k=1}^N \, 
\eps^{\frac{k}{l}} \ U_k ( t,x, \theta) \, , \qquad \qquad \qquad 
\qquad \qquad \ \\
p^\eps_\flat ( t,x, \theta ) \, = \, \p_0 (t,x) + \sum_{k=1}^N \, 
\eps^{\frac{k}{l}} \ P_k ( t,x, \theta) \, .
\end{array} \right. $$
The transition from $ \tilde \u^\eps_\flat $ to $ \u^\eps_\flat $ 
is achieved through the phase shift $ \varphi^\eps_a $

\medskip

$ \tilde U_k ( t,x, \eps^{-1} \, \varphi^\eps_\flat ) \, = \, 
\tilde U_k \bigl( t,x, \eps^{-1} \, \varphi^\eps_g + \varphi_l + 
\sum_{k=l+1}^N \, \eps^{\frac{k}{l} - 1} \ \varphi_k \bigr) \, . $

\medskip

\noindent{Use the Taylor formula in order to absorb the small 
term in the right.} It furnishes the following explicit link 
between the $ (U_k,P_k) $ and the $ (\tilde U_k , \tilde P_k) $
\begin{equation} \label{diction}
\ \left. \begin{array}{l}
U_k \bigl( t,x,\theta - \varphi_l (t,x) \bigr) \, := \, \tilde U_k 
(t,x,\theta) \, + \, \cG^k (\tilde U_1, \cdots, \tilde U_{k-1}) (t,
x,\theta) \, , \ \ \\
P_k \bigl( t,x,\theta - \varphi_l (t,x) \bigr) \, := \, \tilde P_k 
(t,x,\theta) \, + \, \cG^k (\tilde P_1, \cdots, \tilde P_{k-1}) (t,
x,\theta) \, .
\end{array} \right.
\end{equation}
The application $ \cG^k $ can be put in the form

\medskip

$ \cG^k (\tilde U_1, \cdots, \tilde U_{k-1}) \, := \, \sum_{p=1}^{k-1} \, 
\partial_\theta^p \cG^k_p (\tilde U_1, \cdots, \tilde U_{k-p}) \, , 
\qquad k \in \{ 1, \cdots , N \} \, .  $

\medskip

\noindent{The terms $ \cG^k_p $ are given by} 

\medskip

$ \cG^k_p (\tilde U_1, \cdots, \tilde U_{k-p}) \, := \, 
\frac{1}{p \, !} \ \sum_{\alpha \in \cJ^k_p} \ \varphi_{
l+1+\alpha_1} \times \cdots \times \varphi_{l+1+\alpha_p} 
\ \tilde U_{\alpha_{p+1}} \, , $

\medskip

\noindent{where the sum is taken over the set}

\medskip

$ \cJ^k_p \, := \, \bigl \lbrace \, \alpha = (\alpha_1, 
\cdots,\alpha_p, \alpha_{p+1}) \in \NN^{p+1} \, ; $

\smallskip

\hfill $ 0 \leq \alpha_j \leq N-l-1 \, , \qquad \forall \, 
j \in \{1, \cdots, p \} \, , \qquad \qquad \quad \ \ \, $

\smallskip

\hfill $ 1 \leq \alpha_{p+1} \leq k-p \, , \qquad \alpha_1 
+ \cdots + \alpha_p + \alpha_{p+1} = k-p \, \bigr \rbrace 
\, . $

\medskip

\noindent{The relation (\ref{diction}) and the definition of 
$ \cG^k $ imply that}

\medskip

$ \bar U_k (t,x) \, = \, \langle \tilde U_k \rangle (t,x) \, , 
\qquad \forall \, k \in \{ 1, \cdots , N \} \, , \qquad 
\forall \, t \in [0,T] \, . $

\medskip

\noindent{Therefore, prescribing the initial data for the 
$ \bar U_k $ or the $ \langle \tilde U_k \rangle $ amounts 
to the same thing.} The condition (\ref{inipreli}) yields

\medskip

$ \cG^k_p (\tilde U_1, \cdots, \tilde U_{k-p}) (0,x,\theta) 
= 0 \, , \qquad \forall \, k \in \{ 1, \cdots , N \} \, . $

\medskip

\noindent{Since $ \varphi_l (0,\cdot) \equiv 0 $, we have}

\medskip

$ \Pi_0 (0,x) \, U_k^* (0,x,\theta) \, = \, \Pi_0 (0,x) \, 
\tilde U_k^* (0,x,\theta) \, , \qquad \forall \, k \in \{ 1, 
\cdots , N \} \, . $

\medskip

\noindent{It is clearly equivalent to specify the initial 
data for the $ \Pi_0 \, U_k^* $ or the $ \Pi_0 \, \tilde 
U_k^* $.

\bigskip

\noindent{$ \bullet $ {\bf The divergence free relation in 
the variables $ (t,x) $.}} Consider the application 

\medskip

$ \Div \, : \, H^\infty \, \longrightarrow \, \text{\rm Im} \, 
(\Div) \, \subset \, \bigl \lbrace \, \g \in H^\infty \, ; \ \hat \g
(0) = 0 \, \bigr \rbrace \, . $ 

\medskip

\noindent{We can select some special right inverse.}

\begin{lem} \label{invpartbar} There is a linear operator
$ \text{\rm ridiv} \, : \, \text{\rm Im} \, (\Div) \longrightarrow 
H^\infty $ with
\begin{equation} \label{exainv}
\Div \circ \text{\rm ridiv} \ \g \, = \, \g \, , \qquad \forall \, 
\g \in \text{\rm Im} \, (\Div) \, . \qquad \qquad \qquad \quad \ 
\end{equation}
For all $ \iota > 0 $ and for all $ m \in \NN $, there is a 
constant $ C_m^\iota > 0 $ such that
\begin{equation} \label{estiinv}
\parallel \text{\rm ridiv} \ \g \parallel_{H^m} \, \leq \, C_m \ 
\parallel \g \parallel_{H^{m+1+ \frac{d}{2} + \iota}} \, , \qquad 
\forall \, \g \in \text{\rm Im} \, (\Div) \, . \qquad
\end{equation}

\end{lem}

\smallskip

\noindent{\em \underline{Proo}f \underline{o}f \underline{the} 
\underline{Lemma} \underline{\ref{invpartbar}}.} Introduce a 
cut-off function $ \psi \in C^\infty (\RR^d) $ such that

\medskip

$ \bigl \{ \, \xi \, ; \ \psi (\xi) \not = 0 \, \bigr \} \, 
\subset \, B (0,2] \, , \qquad \bigl \{ \, \xi \, ; \ \psi 
(\xi) = 1 \, \bigr \} \, \supset \, B (0,1] \, . $

\medskip

\noindent{For $ g \in \text{\rm Im} \, (\Div) $, take the 
explicit formula}

\medskip

$ \text{ridiv} \, (\g) : = \cF^{-1} \, \bigl( \, \int_0^1 \, 
\nabla_\xi  (\psi \, \hat \g) (r \, \xi) \ dr \, + \, \vert 
\xi \vert^{-2} \ (1-\psi)(\xi) \ \hat \g (\xi) \times \xi \, 
\bigr) \, . $

\medskip
 
\noindent{Since $ \hat \g (0) = 0 $, the relation (\ref{exainv}) 
is satisfied.} For $ s > \frac{d}{2} $, the injection $ H^s 
(\RR^d) \hookrightarrow L^\infty (\RR^d) $ is continuous. It
leads to (\ref{estiinv}). \hfill $ \Diamond $

\bigskip

\noindent{$ \bullet $ {\bf The Leray projector in the variables $ (t,x) $.}} 
Note $ \Pi (\xi) $ the orthogonal projector from $ \RR^d $ onto the plane 

\medskip

$ \xi^\perp \, := \, \lbrace \, u \in \RR^d \, ; \ u \cdot \xi 
= 0 \, \rbrace \, . $

\medskip
 
\noindent{Introduce the closed subspace}

\medskip

$ \text{F} \, := \,  \bigl \lbrace \, \u \in L^2 \, ; 
\ \Div \, \u = 0 \, \bigr \rbrace \, \subset \, L^2 \, . $

\medskip

\noindent{Call $ P $ the orthogonal projector from $ L^2 $ 
onto F.} It corresponds to the Fourier multiplier

\medskip

$ P \, \u \, = \, \Pi (D_x) \, \u \, := \, ( 2 \, \pi)^{- 
\frac{d}{2}} \ \int_{\RR^d} \, e^{i \, x \cdot \xi} \ \Pi 
(\xi) \, \hat \u (\xi) \ d \xi \, . $

\medskip

\noindent{The application $ P $ is the Leray projector onto the 
space of divergence free vector fields.} It is a self-adjoint 
operator such that

\medskip

$ \ker \, \Div \, = \, \text{Im} \, P \, , \qquad \text{Im} \, 
\nabla \, = \, \bigl( \ker \, (\Div) \bigr)^\perp \, = \, \ker \,  
P \, . $

\medskip

\noindent{Consider the Cauchy problem}

\medskip

$ \partial_t \u + \nabla \p = \f \, , \qquad \Div \, \u = 0 \, , 
\qquad \u (0, \cdot) = \h $

\medskip

\noindent{with data $ \f \in L^2_T $ and $ \h \in L^2 $.} 
It leads to the equivalent conditions

\medskip

$ \partial_t \u = P \, \f \, , \qquad \u (0, \cdot) 
= P \, \h \, , \qquad \nabla \p = (\id - P) \, \f \, . $

\bigskip

\noindent{Now we come back to the proof of Theorem \ref{appBKW}.} 
It remains to absorb the term $ \tilde \g^\eps_\flat \in 
\text{Im} \, (\Div) $. To this end, take $ \iota = \frac{1}
{2 \, l} $. Define $ \u^\eps_\flat $ and $ \p^\eps_\flat $ 
as in (\ref{BKWdel}) with the $ U_k $ and $ P_k $ of 
(\ref{diction}). Introduce 

\medskip

$ \cc \u^\eps_\flat \, := \, \r \u^\eps_\flat - \text{ridiv} \, 
\tilde \g^\eps_\flat = \bigcirc (\eps^{\frac{N}{l}-2-\frac{d}{2}}) \, , 
\qquad \cc \p^\eps_\flat \, := \, \r \p^\eps_\flat = \bigcirc 
(\eps^{\frac{N+1}{l}}) \, . $

\medskip

\noindent{After substitution in $ (\cE) $, we lose again a 
power of $ \eps $.} We find

\medskip

$ \f^\eps_\flat = \tilde \f^\eps_\flat - (\text{ridiv} \, 
\tilde \g^\eps_\flat \cdot \nabla) \, \tilde \u^\eps_\flat 
- (\tilde \u^\eps_\flat \cdot \nabla) \, \text{ridiv} \, 
\tilde \g^\eps_\flat $

\smallskip
 
$ \qquad \quad \ - \, \part_t \text{ridiv} \, \tilde \g^\eps_\flat
+ (\text{ridiv} \, \tilde \g^\eps_\flat \cdot \nabla) \, 
\text{ridiv} \, \tilde \g^\eps_\flat \, = \, \bigcirc(\eps^{
\frac{N}{l}-3-\frac{d}{2}}) \, . $

\medskip

\noindent{The Theorem \ref{appBKW} looks like classical statements in one phase 
non linear geometric optics except that the phase $ \varphi^\eps_g $ does depend 
on $ \eps $.} In the next chapter, we examine the part of the $ \varphi_k $ 
which make up $ \varphi^\eps_g $ and $ \varphi^\eps_a $.

\section{The cascade of phases.} {\it Turbulence} and {\it intermittency} are topics
which represent extremely different points of view. Two approaches compete:

\medskip

\noindent{a) }The deterministic approach which study the time evolution
of flows arising in fluid mechanics \cite{B}-\cite{C}-\cite{D}-\cite{E}-\cite{MPP}.

\medskip

\noindent{b) }The statistical approach in which the velocity of the
fluid is a random variable \cite{FMRT}-\cite{L}.

\medskip

\noindent{The Theorem \ref{appBKW} is mainly connected with a)}. It brings
various informations related to the propagation of quasi-singularities. These
aspects are detailed at first. Then we briefly explain b) and we draw 
(in the setting of the Theorem \ref{appBKW}) a phenomenological 
comparison between a) and b).

\smallskip

\subsection{Microstructures.}  The result \ref{appBKW} is concerned
with the convection of microstructures. It is linked with the multiple 
scale approach of \cite{MPP} and \cite{C}. In \cite{MPP} the authors 
look for BKW solutions $  \u^\eps_a $ in the form

\medskip

$ \u^\eps_a (t,x) \, = \, \u_0 (t,x) \, + \, U_0^* \bigl( t,x, 
\eps^{-1} \, t , \eps^{-1} \, \vec \varphi_0 (t,x) \bigr) \, + \, 
\bigcirc (\eps) \, . $ 

\medskip

\noindent{In the more recent paper \cite{C}, the selected expansion is}

\medskip

$ \mu^\eps_a (t,x) \, = \, \u_0 (t,x) \, + \, \eps^{\frac{1}{3}} \ 
U_1 \bigl( t,x, \eps^{-\frac{2}{3}} \, t, \eps^{-1} \, \vec \varphi_0 
(t,x) \bigr) \, + \, \bigcirc (\eps^{\frac{2}{3}}) \, . $ 

\medskip

\noindent{Both articles  \cite{C} and \cite{MPP} use homogenization 
techniques.} They perform computations involving expressions as 
$ \u^\eps_a $ or $ \mu^\eps_a $. Simplifications (supported by 
engineering experiments) are made in order to get effective 
equations for the evolution of $ (\u_0 , U_0^*) $ or $ (\u_0 ,
U_1 ) $. 

\smallskip

\noindent{Consider the simple case of one phase expansions (that 
is when $ \vec \varphi_0 \equiv \varphi_0 $ is a scalar valued 
function).} Reasons why a complete mathematical analysis based on 
$ \u^\eps_a $ or $ \mu^\eps_a $ is not available can be drawn from 
the Theorem \ref{appBKW}. For instance, look at $ \mu^\eps_a $.
When $ l = 3 $, the oscillation $ \mu^\eps_a $ involves the same 
scales as $ \u^\eps_{(3,N)} $ since

\medskip

$ \eps^{-1} \, \varphi^\eps_g (t,x) \, = \, \eps^{-1} \ \varphi_0 
(t,x) + \eps^{- \frac{2}{3}} \ \varphi_1 (t,x) + \eps^{- \frac{1}
{3}} \  \varphi_2 (t,x) \, . $

\medskip

\noindent{Now the analogy stops here since in general $ \varphi_1 
(t,x) \not \equiv t $ and $ \varphi_2 (t,x) \not \equiv 0 $.} These 
are geometrical obstructions which prevent to describe the propagation 
by way of $ \mu^\eps_a $. The asymptotic expansion $ \mu^\eps_a $ 
is not suitable.

\medskip

\noindent{Analogous arguments concerning $ \u^\eps_a $ will be 
presented in the paragraph 3.5.}

\smallskip

\subsection{The geometrical phase.} Let us examine 
more carefully how the expression $ \varphi^\eps_g $ 
is built. Because of the condition (\ref{nonstaou}), 
for $ \eps $ small enough, it is still not stationary
\begin{equation} \label{nondegeat}
\ \exists \, \eps_0 > 0 \, ; \qquad \nabla \varphi^\eps_g 
(t,x) \not = 0 \, , \qquad \forall \, (\eps,t,x) \in \, 
]0,\eps_0] \times [0,T] \times \RR^d \, . \ \
\end{equation}
In fact, the function $ \varphi^\eps_g $ comes from the 
approximate eiconal equation

\medskip

$ \partial_t \varphi^\eps_g + ( \bar u^\eps_\flat \cdot 
\nabla) \varphi^\eps_g = \bigcirc ( \eps ) $

\medskip

\noindent{which is equivalent to}

\medskip

$ \partial_t \varphi_k + \u_0 \cdot \nabla \varphi_k + \sum_{j=0}^{k-1}
\, \bar U_{k-j} \cdot \nabla \varphi_j = 0 \, , \qquad \forall \, k \in \{ 1,
\cdots, l-1 \} \, .  $

\medskip

\noindent{The family $ \{ \u^\eps_\flat (t,x) \}_{\eps \in \, ]0,1]} $
has an $ \eps - \, $stratified regularity \cite{G2} with respect to the
phase $ \varphi^\eps_g $.} This is a {\it geometrical} information.

\smallskip

\subsection{Closure problems.} We have explained why appealing 
only to $ \varphi_0 $ is not sufficient. It turns out that BKW 
computations relying only on the geometrical phase $ \varphi^\eps_g $ 
come also to nothing. This is a subtle aspect when proving the 
Theorem \ref{appBKW}. We lay now stress on it.

\medskip

\noindent{For all $ N \in \NN_* $, the application $ \cG $ 
defined below is one to one}
$$ \ \left. \begin{array}{rcl}
\cG \ \, : \, (\cW^\infty_T)^N & \longrightarrow & (\cW^\infty_T)^N \\
\\
\left(\begin{array}{c}
\tilde U_1 \\
\tilde U_2 \\
\vdots \\
\tilde U_N
\end{array} \right)(t,x,\theta) & \longmapsto & 
\left(\begin{array}{c}
\tilde U_1 \\
\tilde U_2 + \cG^1 (\tilde U_1) \\
\vdots \\
\tilde U_N + \cG^N (\tilde U_1, \cdots, \tilde U_{N-1}) 
\end{array} \right)(t,x,\theta + \varphi_l (t,x) ) \, .
\end{array} \right. $$
Once the $ U_j $ or the $ \tilde U_j $ are known, it is entirely 
equivalent to use $ \u^\eps_\flat $ or $ \tilde \u^\eps_\flat $. 
Before the $ U_j $ or the $ \tilde U_j $ have been identified, in 
particular when performing the BKW calculus, it is deeply different to 
employ $ \u^\eps_\flat $ or $ \tilde \u^\eps_\flat $. Indeed, there is 
a unique choice of the $ \varphi_k $ with $ l \leq k \leq N $, which 
imposes a {\it specific hierarchy} between the profiles $ \tilde 
U_k $, which makes possible the {\it triangulation} of the equations 
obtained by the formal computations.

\medskip

\noindent{In the subsection 2.3, we will perform the BKW analysis
with the profiles $ \tilde U_k $.} It yields a sequence of equations
\begin{equation} \label{cloeetilde}
\tilde X^k ( \tilde U_1, \cdots, \tilde U_{k+l} ) \, = \, 0 \, , 
\qquad 1 \leq k \leq N \, . \qquad \qquad \quad \
\end{equation}
As usual in non linear geometric optics, this can be rewritten 
in order to find a sequence of well-posed equations
\begin{equation} \label{cloeq}
\dot X^k ( \dot U_k ) \, = \, \cF ( \dot U_1, \cdots , \dot 
U_{k-1} ) \, , \qquad 1 \leq k \leq N \, , \qquad \qquad
\end{equation}
where the $ \dot U_k $ are made of pieces of the $ \tilde
U_j $. Of course, the equation (\ref{cloeq}) can be interpreted 
in terms of the $ \tilde U_j $ and then in terms of the $ U_j $. 
In this second step, it requires to implement the phase 
shift $ \varphi_l $ and the transformations $ \cG^j_p $ 
with $ 1 \leq j \leq k-1 $ and $ 1 \leq p \leq j $. Now, 
the BKW analysis reveals that $ \varphi_l $ or the various 
coefficients $ \varphi_i $ which appear in the definition 
of such $ \cG^j_p $ do not depend only on $ (\dot U_1, 
\cdots, \dot U_k) $ but also on some $ \dot U_i $ with 
$ i > k $. The resulting system is therefore underdetermined.
Computations involving the functions $ U_j $ lead to a 
sequence of equations which are not closed.

\medskip

\noindent{The insertion of the phases $ \varphi_k $ with $ 1 
\leq k \leq N $ is an elegant way to introduce $ \cG $.} The 
change of variables $ \cG $, though it is a function of $ (U_1, 
\cdots, U_N) $, is needed to progress. It allows to get round 
{\it closure problems}.

\subsection{Compensated compactness.} Dissipation terms can be 
incorporated in the discussion. In the variables $ (t,x) $, the 
addition of some viscosity $ \kappa $ is compatible with the 
propagation of oscillations if for instance $ \kappa = \nu \, 
\eps^2 $. There are approximate solutions $ (\u^\eps_\flat ,
\p^\eps_\flat ) $ of the Navier-Stokes equations. They satisfy
(\ref{BKWdel}) and 
\medskip

$ \part_t \u^\eps_\flat + ( \u^\eps_\flat \cdot \nabla) \u^\eps_\flat 
+ \nabla \p^\eps_\flat = \nu \ \eps^2 \ \Delta_x \u^\eps_\flat + 
\f^\eps_\flat \, , \qquad \Div \ \u^\eps_\flat = 0 \, , $

\medskip

\noindent{with $ \f^\eps_\flat = \bigcirc(\eps^\infty) $.}
When $ \nu > 0 $, Leray's theorem provides with global weak 
solutions $ (\u^\eps,\p^\eps) (t,x) $ of the following 
Cauchy problem}
$$ \left \{ \begin{array} {l}
\part_t \u^\eps + ( \u^\eps \cdot \nabla) \u^\eps + \nabla 
\p^\eps = \nu \ \eps^2 \ \Delta_x \u^\eps \, , \qquad \Div \ 
\u^\eps = 0 \, , \qquad \qquad \\
\u^\eps (0,\cdot) \equiv \u^\eps_\flat (0,\cdot) \, . 
\end{array} \right. $$
Suppose now that $ \u_0 \equiv 0 $. Then, we have 
also the uniform controls
\begin{equation} \label{compcomp}
\, \left. \begin{array}{l}
\sup \ \bigl \lbrace \, \parallel \eps^{-\frac{1}{l}} \,
\u^\eps \parallel_{L^2_T} \, ; \ \eps \in \, ]0,1] \, 
\bigr \rbrace \, \leq \, C \, < \, \infty \, , \\
\sup \ \bigl \lbrace \, \nu \ \eps^2 \ \int_0^T \, 
\parallel \eps^{-\frac{1}{l}} \, \u^\eps (t,\cdot) 
\parallel^2_{H^1 (\RR^d) } \ dt \,  ; \ \eps \in \, 
]0,1] \, \bigr \rbrace \, \leq \, C \, < \, \infty \, . \
\end{array} \right.
\end{equation}
Arguments issued from the theory of compensated compactness 
\cite{Ge} can be employed to study the sequence $ \{ \eps^{
-\frac{1}{l}} \, \u^\eps \}_\eps $. In the spirit of \cite{D} 
or \cite{E}, we can try to exploit the informations contained 
in (\ref{compcomp}) and the equation on $ \u^\eps $ in order 
to describe the asymptotic behaviour when $ \eps $ goes to 
zero of the functions $ \eps^{- \frac{1}{l}} \, \u^\eps $. 
However this approach seems to be not applicable here.

\medskip

\noindent{Indeed, {\it obvious} instabilities occur.} 
The related mechanisms, which induce the non linear instability 
of Euler equations, are detailed in the paragraph 5.1. Below, we 
just give an intuitive idea of what can happen. Use the representation 
$ \tilde \u^\eps_\flat $ involving the phase $ \varphi^\eps_\flat $. 
The determination of the intermediate term $ \varphi_l $ requires to 
identify $ \langle \tilde U_l \rangle $ and $ \tilde U^*_{l-1} $. This 
is a consequence of the equations (\ref{phik}) and (\ref{moyj+1}). 

\medskip

\noindent{In view of the formula (\ref{diction}), when $ \varphi_l $ 
is modified by an amount of $ \delta \varphi_l $, the quantity 
$ U_1 (t,x,\theta) $ undergoes a perturbation of the same order 
$ \delta \varphi_l $.} When dealing with quasi-singularities, 
some quantities with $ \eps $ in factor (like $ \langle \tilde 
U_l \rangle $) or with $ \eps^{1 - \frac{1}{l}} $ in factor 
(like $ \tilde U^*_{l-1} $) can control informations of size 
$ \eps^{\frac{1}{l}} $. This fact is expressed by the 
following rules of transformation
\begin{equation} \label{rulesoftr}
\left. \begin{array}{lcl}
\langle \tilde U_l \rangle \ / \ \langle \tilde U_l \rangle \, 
+ \, \delta \langle \tilde U_l \rangle \quad & 
\Longrightarrow & \quad \u^\eps_\flat \ / \ \u^\eps_\flat + 
\bigcirc ( \eps^{\frac{1}{l}} ) \ \delta \langle \tilde U_l 
\rangle \, , \\
\tilde U^*_{l-1} \ / \ \tilde U^*_{l-1} \, + \, \delta 
\tilde U^*_{l-1} \quad & \Longrightarrow & \quad \u^\eps_\flat \ 
/ \ \u^\eps_\flat + \bigcirc ( \eps^{\frac{1}{l}} ) \ \delta 
\tilde U^*_l \, . \
\end{array} \right. 
\end{equation}
Now reverse the preceding reasoning. To describe features in the
principal oscillating term $ \eps^{\frac{1}{l}} \ U^*_1 
\bigl(t,x,\eps^{-1} \, \varphi^\eps_g (t,x) \bigr) $, we 
must identify $ \varphi_l $ which means to obtain $ \langle 
\tilde U_l \rangle $ and $ \tilde U^*_{l-1} $. In other words, 
we need to know quantities which have respectively $ \eps $ and 
$ \eps^{1- \frac{1}{l}} $ in factor. When $ l \geq 2 $ such 
informations are clearly not reachable by rough controls as 
(\ref{compcomp}). 

\medskip

\noindent{This discussion indicates that the study of turbulent 
regimes requires to combine at least geometrical aspects, multiphase 
analysis and high order expansions.} The tools of non linear geometric 
optics seem to be appropriate. Some attempts in this direction have 
already been made.

\smallskip

\subsection {Non linear geometric optics.} We make in this paragraph 3.5
several comments about non linear geometric optics. They concern both old
\cite{G}-\cite{G2}-\cite{Se} and recent \cite{Che}-\cite{CGM}-\cite{CGM1} 
results which all are devoted to one phase expansions of the type
\begin{equation} \label{ogmonop}
\u^\eps_\natural (t,x) \, := \, \u_0 (t,x) \,  + \, \hbox{$ \sum_{k=1}^\infty $}
\, \eps^{\frac{k}{l}} \ U_k \bigl( t,x, \eps^{-1} \, \varphi_0 (t,x) \bigr) 
\, . \quad \
\end{equation}
When $ l=1 $, one is faced with {\it weakly non linear geometric optics}. 
The asymptotic behavior and the stability of $ \u^\eps_\natural $ are
well understood. In fact a complete theory has been achieved in the
general framework of multidimensional systems of conservation laws 
(see \cite{G}-\cite{G2} and the related references). Because of the
formation of shocks, the life span of exact solutions close to 
$ \u^\eps_\natural $ does not go beyond $ T \simeq 1 $. 

\medskip

\noindent{When $ l=2 $, expressions as $ \u^\eps_\natural $ are called 
{\it strong} oscillations.} The construction of such BKW solutions can 
be undertaken only if the system of conservation laws has a special 
structure. {\it Transparency} conditions are needed to progress. They
can be deduced from the presence of a linearly degenerate field
\cite{CGM}. In the hyperbolic situation the family $ \{ \u^\eps_\natural 
\}_{\eps \in \, ]0,1]} $ is unstable \cite{CGM} on the interval $ [0,T] $.
It becomes stable on condition that a small viscosity is incorporated
\cite{Che}. Applications can be given to describe large-scale 
motions in the atmosphere \cite{Che}.

\medskip

\noindent{Compressible Euler equations are the prototype of a non 
linear hyperbolic system having a linearly degenerate field.} After
a finite time, singularities appear. These correspond to the generation
of shocks by compression \cite{Si}. The situation is different in the
incompressible setting. There is no genuine shock and the production
of singularities poses a much more subtle problem \cite{BKM}-\cite{CF}
which up to now remains basically open.

\medskip

\noindent{Incompressible fluid equations lie at an extreme end in the 
sense that they are the most {\it degenerate} (or the most {\it linear})
equations which have just been mentioned. Following the approach of 
\cite{JMR1} related to {\it transparency}, repeating the reasoning 
which goes from \cite{G}-\cite{G2} to \cite{Che}-\cite{CGM}, one 
expects to go further than $ l = 2 $ when dealing with $ (\cE) $. 
Now, this is precisely what says the Theorem \ref{appBKW} since 
it allows to reach any $ l \in \NN_* \, !$

\medskip

\noindent{To tackle the limit case $ l= \infty $, one is tempted to look
at asymptotic expansions of the form
\begin{equation} \label{ogmonopl}
\ \u^\eps_\infty (t,x) \, := \, \hbox{$ \sum_{k=0}^\infty $} \, \eps^k \ 
U_k \bigl( t,x, \eps^{-1} \, \varphi_0 (t,x) \bigr) \, , \qquad
\partial_\theta U_0^* \not \equiv 0 \, . \qquad \quad  
\end{equation}
The oscillations contained in $ \u^\eps_\infty $ have a {\it large}
amplitude. Modulation equations for $ U_0 $ are proposed in \cite{Se}. 
However these transport equations are not hyperbolic so that they are 
ill posed (in the sense of Hadamard) with respect to the initial value 
problem. It confirms that a BKW construction based on (\ref{ogmonopl}) 
is not relevant\footnote{The singularities are carried here by the velocity 
field. The discussion is very different when the oscillations are 
polarized on the entropy \cite{CGM1}.}. 

\medskip

\noindent{The contribution \cite{Se} does not explain why the expansion
(\ref{ogmonopl}) is not the good one.} We come back below to this point.
At first sight the Theorem \ref{appBKW} does not include large 
amplitude waves since $ \u^\eps_\flat - \u_0 = \bigcirc ( \eps^{\frac
{1}{l}} ) \ll \bigcirc(1) $. A change of variables leads to recant 
this impression. Suppose that $ \u_0 \equiv 0 $ and $ \part_\theta 
U_1^* \not \equiv 0 $. Then define
 
\medskip

$ \dot \u^\eps_\flat (t,x) \, := \, \eps^{-\frac{1}{l}} \ \u^\eps_\flat 
(\eps^{-\frac{1}{l}} \, t,x) \, , \qquad \dot \p^\eps_\flat (t,x) \, :=
\, \eps^{-\frac{2}{l}} \ \p^\eps_\flat (\eps^{-\frac{1}{l}} \, t,x) \, . $

\medskip

\noindent{Observe that the structure of $ \dot \u^\eps_\flat $ and 
$ \dot \p^\eps_\flat $ is very different from the one in (\ref{ogmonopl}) 
since we have
$$ \quad \left. \begin{array}{l}
\dot \u^\eps_\flat (t,x) = \sum_{k=1}^\infty \, \eps^{ 
\frac{k-1}{l}} \ U_k \bigl( \eps^{-\frac{1}{l}} \, t,x, \eps^{-1} \, 
\varphi^\eps_g (\eps^{-\frac{1}{l}} \, t,x) \bigr) + 
\eps^{-\frac{1}{l}} \ \cc \u^\eps_\flat (\eps^{-\frac{1}{l}} \, t,x) 
 \, , \\
\dot \p^\eps_\flat (t,x) = \sum_{k=1}^\infty \, \eps^{
\frac{k-2}{l}} \ P_k \bigl( \eps^{-\frac{1}{l}} \, t,x, \eps^{-1} \, 
\varphi^\eps_g (\eps^{-\frac{1}{l}} \, t,x) \bigr) + 
\eps^{-\frac{2}{l}} \ \cc \p^\eps_\flat (\eps^{-\frac{1}{l}} \, t,x) \, .
\end{array} \right. $$
The functions  $ \dot \u^\eps_\flat $ and $ \dot \p^\eps_\flat $
satisfy

\medskip

$ \part_t \dot \u^\eps_\flat + ( \dot \u^\eps_\flat \cdot \nabla) 
\dot \u^\eps_\flat + \nabla \dot \p^\eps_\flat = \dot \f^\eps_\flat \, , 
\quad \ \Div \, \dot \u^\eps_\flat = 0 \, , \quad \ \dot \f^\eps_\flat 
(t,x) = \eps^{-\frac{2}{l}} \ \f^\eps_\flat (\eps^{-\frac{1}{l}} \, t,x) 
\, . $

\medskip

\noindent{The functions $ \dot \u^\eps_\flat $ are oscillations 
of the order $ 1 $.} They are approximate solutions of $ (\cE) $ 
on the {\it small} interval $ [0, \eps^{\frac{1}{l}} \, T] $.
Indeed, for all $ m \in \NN $, the family $ \{ \dot \f^\eps_\flat 
\}_\eps $ is subjected to the uniform majoration

\medskip

$ \sup_{\eps \in \, ]0,\eps_0]} \quad \eps^{- \frac{N}{l} +
\frac{2}{l} + 3 + m} \ \parallel \dot \f^\eps_\flat \parallel_{
\cW^m_{\! \text{\tiny $ \eps^{(1/l)} \, T $}} } \ < \, \infty \, . \ $

\medskip

\noindent{If moreover $ N = + \infty $ and} 
\begin{equation} \label{speini}
\left. \begin{array}{l}
\varphi_1 (0,\cdot) \, \equiv \, \cdots  \, \equiv \, \varphi_{l-1}
(0,\cdot) \, \equiv \, 0 \, , \\
U_{k+1} (0,\cdot) \, \equiv \, 0 \, , \qquad \forall \, k \in \NN 
\setminus ( l \, \NN) \, , \qquad \qquad \qquad \qquad \qquad \
\end{array} \right.
\end{equation}
the trace $ \dot \u^\eps_\flat (0,\cdot) $ has the form

\medskip

$ \dot \u^\eps_\flat (0,x) \, = \ \sum_{k=0}^\infty \, \eps^k \  
U_{1 + l \, k} \bigl(0,x, \eps^{-1} \, \varphi_0 (0,x) \bigr) \, , 
\qquad \part_\theta U_1^* \, \not \equiv \, 0 \, .$

\medskip

\noindent{At the time $ t=0$, we recover (\ref{ogmonopl}). Now 
the construction underlying the Theorem \ref{appBKW} reveals that 
in general}
\begin{equation} \label{phasenonn}
\varphi_k(t,\cdot) \, \not \equiv \, 0 \, , \qquad \forall \, t \in \, ]0,T] \, , \qquad 
\forall \, k \in \{ 2, \cdots, l-1 \} \, . \quad 
\end{equation}
The functions $ \varphi_j $ with $ j \in \{2,\cdots, l-1 \} $ are 
not present when $ t = 0 $. But the description of $ \dot \u^\eps_\flat 
(t,\cdot) $ on the interval $ [0, \eps^{1 - \frac{k}{l}} \, T] $ with $ 
k \in \{2,\cdots, l-1 \} $ requires the introduction of the {\it phase 
shifts} $ \varphi_j $ for $ j \in \{2,\cdots,k\} $. More generally, the 
description of $ \dot \u^\eps_\flat (t,\cdot) $ on the whole interval 
$ [0, T] $ needs the introduction of an {\it infinite cascade} of phases 
$ \{ \varphi_j \}_{j \in \NN_*} $.

\medskip

\noindent{Such a phenomenon does not occur when constructing large 
amplitude oscillations for systems of conservation laws in one space 
dimension \cite{CG}-\cite{E2}.} It is specific to the multidimensional 
framework. It explains why the classical approach of \cite{Se} fails.

\medskip

\noindent{It seems that the creation of the $ \varphi_j $ is due to 
mechanisms which have not already been studied.} It is not linked with 
{\it resonances}. It is related neither to {\it dispersive} nor to 
{\it diffractive} effects.

\bigskip
\vskip -1mm

\noindent{{\it Remark 3.5.1 (about $ \varphi_1 $):}} The term $ \varphi_1 $ 
does not appear if $ \varphi_1 (0,\cdot) \equiv 0 $ and $ \bar U_1
(0,\cdot) \equiv 0 $. When these two conditions are not verified, the 
phase shift $ \varphi_1 $ can be absorbed by the technical trick 
exposed in \cite{CGM}. Just replace $ \u_0 (0,\cdot) $ by $ \u_0 
(0,\cdot) + \delta \, \bar U_1 (0,\cdot) $. Perform the BKW calculus 
with a fixed $ \delta > 0 $. Then choose $ \delta = \eps $.  \hfill 
$ \triangle $

\bigskip

\noindent{{\it Remark 3.5.2 (about $ \varphi_2 $):}} In general, we
have $ \varphi_2 \not \equiv 0 $ even if

\medskip

$ \varphi_1 (0,\cdot) \equiv \varphi_2 (0,\cdot) \equiv 0 \, , \qquad 
\bar U_1 (0,\cdot) \equiv \bar U_2 (0,\cdot) \equiv 0 \, . $

\medskip

\noindent{Indeed the time evolution of $ \bar U_2 $ is governed 
by (\ref{moy2}).} It involves the source term $ \Div \, \langle 
U^*_1 \otimes U^*_1 \rangle $ which is able to awake $ \bar U_2 $. This 
influence can then be transmitted to $ \varphi_2 $ through the 
transport equation
\begin{equation} \label{trans2}
\part_t \varphi_2 + (\u_0 \cdot \nabla) \varphi_2 + (\bar U_1 
\cdot \nabla) \varphi_1 + (\bar U_2 \cdot \nabla) \varphi_0
\, = \, 0 \, . \qquad \qquad
\end{equation}
Likewise, the other terms $ \varphi_3 $, $\cdots $, 
$ \varphi_{l-1} $ are in general non trivial even if 

\medskip

$ \varphi_1 (0,\cdot) \equiv \cdots \equiv \varphi_{l-1} (0,\cdot) 
\equiv 0 \, , \qquad \bar U_1 (0,\cdot) \equiv \cdots \equiv 
\bar U_{l-1} (0,\cdot) \equiv 0 \, . $

\medskip

\noindent{There is no more trick which allows to get rid of 
$  \varphi_2 $, $ \cdots$, $ \varphi_{l-1} $.} \hfill $ \triangle $

\bigskip

\noindent{{\it Remark 3.5.3 (why turbulent flows ?):}} The introduction 
of the phase shifts $  \varphi_k $ with $ 2 \leq k \leq l-1 $ cannot 
be avoided. Therefore the difficulties that we deal with appear from 
$ l = 3 $. When $ l \geq 3 $, the characteristic rate $ e $ of eddy  
dissipation is bigger than one \cite{C}. This is the reason why such 
situations are refered to {\it turbulent regimes}. \hfill $ \triangle $

\bigskip

\noindent{{\it Remark 3.5.4 (about shear layers):}} We have said in 
the introduction that the expression $ \u^\eps_s $ given by formula
(\ref{oscidipmaj}) is of a very special form. Let us explain why.
Change the variable $ t $ into $ \eps^{\frac{1}{l}} \, t $ 
and $ \u^\eps_s $ into $ \dot \u^\eps_s := \eps^{\frac{1}{l}} \, 
\u^\eps_s $. The main phase $ \varphi_0 (t,x) \equiv x_2 $ remains 
the same. Now we are faced with

\medskip

$ \dot \u^\eps_s (t,x) \, := \, {}^t \bigl( \eps^{\frac{1}{l}} \, 
\g ( x_2,\eps^{-1} \, x_2) , 0 , \eps^{\frac{1}{l}} \, \h \bigl( 
x_1 - \eps^{\frac{1}{l}} \, \g (x_2, \eps^{-1} \, x_2) \, t , x_2, 
\eps^{-1} \, x_2 \bigr) \bigr) \, . $

\medskip

\noindent{It is still a solution of Euler equations.} Now it falls
in the framework of the Theorem \ref{appBKW}. The constraints 
on $ \bar U_2 = {}^t (\bar U_2^1, \bar U_2^2, \bar U_2^3 ) $ 
reduce to

\medskip

$ \bar U_2^1 \, \equiv \, \bar U_2^2 \, \equiv \, 0 \, , \qquad 
\part_t \bar U_2^3 \, + \, \langle \g \, \part_1 \h \rangle \,
= \, 0 \, . $

\medskip

\noindent{The contribution $ \bar U_2 $ is non trivial but it is
polarized so that $ \bar U_2 \cdot \nabla \varphi_0 \equiv 0 $.}
Therefore it does not produce the phase shift $ \varphi_2 $. The
same phenomenon occurs concerning $ \varphi_3 $, $ \cdots , $ 
$ \varphi_{l-1} $. These terms are not present. It turns out that 
the expansion $ \u^\eps_s $ involves only the phase $ \varphi_0 
(t,x) \equiv x_2 $. \hfill $ \triangle $

\bigskip

\noindent{The choice for the amplitude of the oscillations is very 
important.} It is strongly related to the scale of time $ T $ under
consideration. The idea is to increase the time of propagation $ T $
to reach the regime where non linear effects appear. Starting
with some large amplitude high frequency waves 

\medskip

$  \u^\eps_\infty (0,x) \, = \, U_0 \bigl( 0,x, \eps^{-1} \, \varphi_0 
(0,x) \bigr) \, + \, \bigcirc (\eps) \, , \qquad \partial_\theta U_0^* 
(0,\cdot) \not \equiv 0 \, , $

\medskip

\noindent{the preceding discussion can be summarized by the following 
diagram:}

\bigskip

$$ \left. \begin{array}{ccccccc}
T \simeq 1 &
-\!\!\!-\!\!\!-\!\!\!-\!\!\!-\!\!\!-\!\!\!-\!\!\!-\!\!\!-\!\!\!-\!\!\!-
& \ & 
-\!\!\!-\!\!\!-\!\!\!-\!\!\!-\!\!\!-\!\!\!-\!\!\!-\!\!\!-
& \ & 
-\!\!\!-\!\!\!-\!\!\!-\!\!\!-\!\!\!-\!\!\!-\!\!\!-\!\!\!-\!\!\!-\!\!\!-\!\!\!-\!\!\!-\!\!\!-\!\!\!-\!\!\!-\!\!\!-\!\!\!-\!\!\!- \\

\ & \left. \begin{array}{c}
\text{\scriptsize infinite cascade} \\
\text{\scriptsize of phases} \\
\text{\scriptsize  $\varphi_0 - (\varphi_1) - \cdots $}
\end{array} \right. & \left. \begin{array}{c}
\mid \\
\mid \\ 
\mid
\end{array} \right. &  \left. \begin{array}{c}
 \\
\text{\scriptsize turbulent} \\ 
\text{\scriptsize flows}
\end{array} \right.

& \left. \begin{array}{c}
\mid \\
\mid \\ 
\mid
\end{array} \right.
 & \left. \begin{array}{c}
\\
\text{\footnotesize incompressible} \\
\text{\footnotesize fluid equations}
\end{array} \right.
  \\

T \simeq \eps^{\frac{1}{3}} \ &
 -\!\!\!-\!\!\!-\!\!\!-\!\!\!-\!\!\!-\!\!\!-\!\!\!-\!\!\!-\!\!\!-\!\!\!-
& \ &
-\!\!\!-\!\!\!-\!\!\!-\!\!\!-\!\!\!-\!\!\!-\!\!\!-\!\!\!-
& \ & 
-\!\!\!-\!\!\!-\!\!\!-\!\!\!-\!\!\!-\!\!\!-\!\!\!-\!\!\!-\!\!\!-\!\!\!-\!\!\!-\!\!\!-\!\!\!-\!\!\!-\!\!\!-\!\!\!-\!\!\!-\!\!\!- \\

\ & \text{\scriptsize  $\varphi_0 - (\varphi_1) - \varphi_2 $} & \left. \begin{array}{c}
\mid \\
\mid
\end{array} \right. & \left. \begin{array}{c}
\text{\footnotesize turbulent} \\
\text{\footnotesize flows}
\end{array} \right. & \left. \begin{array}{c}
\mid \\
\mid
\end{array} \right. & \left. \begin{array}{c}
\text{\footnotesize incompressible} \\
\text{\footnotesize fluid equations}
\end{array} \right. \\

T \simeq \eps^{\frac{1}{2}} \ & -\!\!\!-\!\!\!-\!\!\!-\!\!\!-\!\!\!-\!\!\!-\!\!\!-\!\!\!-\!\!\!-\!\!\!-
& \ &
-\!\!\!-\!\!\!-\!\!\!-\!\!\!-\!\!\!-\!\!\!-\!\!\!-\!\!\!-
& \ & -\!\!\!-\!\!\!-\!\!\!-\!\!\!-\!\!\!-\!\!\!-\!\!\!-\!\!\!-\!\!\!-\!\!\!-\!\!\!-\!\!\!-\!\!\!-\!\!\!-\!\!\!-\!\!\!-\!\!\!-\!\!\!- \\

&\text{\footnotesize $ \varphi_0 - (\varphi_1) $} & \left. \begin{array}{c}
\mid \\
\mid \\ 
\mid
\end{array} \right. & \left. \begin{array}{c}
\text{\small strong}\\
\text{\small oscillations} \\
\cite{Che}-\cite{CGM}
\end{array} \right. & \left. \begin{array}{c}
\mid \\
\mid \\ 
\mid
\end{array} \right. & \left. \begin{array}{c}
\text{\small systems of conservation}\\
\text{\small laws with a linearly} \\
\text{\small degenerate field}
\end{array} \right.
 \\
T \simeq \eps & -\!\!\!-\!\!\!-\!\!\!-\!\!\!-\!\!\!-\!\!\!-\!\!\!-\!\!\!-\!\!\!-\!\!\!-
& \ &
-\!\!\!-\!\!\!-\!\!\!-\!\!\!-\!\!\!-\!\!\!-\!\!\!-\!\!\!-
& \ & -\!\!\!-\!\!\!-\!\!\!-\!\!\!-\!\!\!-\!\!\!-\!\!\!-\!\!\!-\!\!\!-\!\!\!-\!\!\!-\!\!\!-\!\!\!-\!\!\!-\!\!\!-\!\!\!-\!\!\!-\!\!\!- \\

 & \varphi_0 & \left. \begin{array}{c}
\mid \\
\mid \\ 
\mid \\
\mid \\
\mid
\end{array} \right. &  \left. \begin{array}{c}
\text{weakly} \\
\text{non linear} \\
\text{geometric}\\
\text{optics} \\ 
 \cite{G}-\cite{G2}
\end{array} \right.
 & \left. \begin{array}{c}
\mid \\
\mid \\ 
\mid \\
\mid \\
\mid
\end{array} \right. &  \left. \begin{array}{c}
\text{\large systems} \\
\text{\large of} \\
\text{\large conservation} \\ 
\text{\large laws}
\end{array} \right. \\

T=0 & -\!\!\!-\!\!\!-\!\!\!-\!\!\!-\!\!\!-\!\!\!-\!\!\!-\!\!\!-\!\!\!-\!\!\!-
& \ &
-\!\!\!-\!\!\!-\!\!\!-\!\!\!-\!\!\!-\!\!\!-\!\!\!-\!\!\!-
& \ & -\!\!\!-\!\!\!-\!\!\!-\!\!\!-\!\!\!-\!\!\!-\!\!\!-\!\!\!-\!\!\!-\!\!\!-\!\!\!-\!\!\!-\!\!\!-\!\!\!-\!\!\!-\!\!\!-\!\!\!-\!\!\!- \\

& \text{\large \bf phases} & \mid & \text{\large \bf regimes} & \mid & \text{\large \bf equations} \\

\end{array} \right. \quad $$

\bigskip

\noindent{This picture allows to understand the position of the
actual paper in comparison with previous results.}

\smallskip

\subsection{The statistical approach.} It deals mainly with
{\it quantitative} informations obtained at the level of expressions, 
say $ \mu (x) $, which in general do not depend on the time $ t $. 
The introduction of $ \mu $ can be achieved by looking at {\it 
stationary statistical} solutions \cite{FMRT} of the Navier-Stokes  
equations that is

\medskip

$ \mu(x) \, \equiv \, \lim_{\, T \, \longrightarrow \, \infty} \quad
\frac{1}{T} \ \int_0^T \, \u(t,x) \ dt $

\medskip

\noindent{or in conjunction with the {\it ensemble average operator} 
(\cite{L}-V-6) marked by the brackets $ < \cdot > $.} We will follow 
this second option. The description below is extracted from the book 
of M. Lesieur \cite{L} (chapters V and VI). We work with $ d = 3 $. 
Interesting quantities are the mean kinetic energy

\medskip

$ \frac{1}{2} \ < \mu(x)^2 > \ \sim \, \int_{\RR^3} \, \vert \mu(x) \vert^2 \ dx \, , $

\medskip

\noindent{the enstrophy (that is the space integral of the square norm of the vorticity)}

\medskip

$ \frac{1}{2} \ < \omega (x)^2 > \ \sim \, \int_{\RR^3} \, \vert \omega (x) 
\vert^2 \ dx \, , \qquad \omega (x) := \nabla \wedge \mu (x) $

\medskip

\noindent{and the rate of dissipation $ e \, \sim \, \kappa \ < 
\omega (x)^2 > $.} In the setting of {\it isotropic} turbulence, 
these quantities can be expressed in terms of a scalar function 
$ k \longmapsto E(k) $. The real number $ E(k) $ represents the 
density of kinetic energy at wave number $ k $ (or the kinetic 
energy in Fourier space integrated on a sphere of radius $ k $). 
The relations are the following

\medskip

\noindent{\cite{L}}-V-10-4$ \qquad \ \frac{1}{2} \ < \mu(x)^2 > \ = \, \int_0^{+ \infty} \, 
E(k) \ dk \, . $ 

\medskip

\noindent{\cite{L}}-V-10-15$ \qquad \! \frac{1}{2} \ < \omega(x)^2 > \ = \, \int_0^{+ \infty} \, 
k^2 \ E(k) \ dk \, . $ 

\medskip

\noindent{\cite{L}}-VI-3-15$ \qquad e \, = \,  2 \ \kappa \ \int_0^{+ \infty} \, k^2 \  
E(k) \ dk \, . $ 

\medskip

\noindent{Kolmogorov's theory assumes that}

\medskip

\noindent{\cite{L}}-VI-4-1 $ \qquad \exists \ c > 0 \, ; \qquad E(k) \,= \, c \ e^{2/3} \ 
k^{-5/3} \, , \qquad \forall \, k \in [k_i,k_d] \, . $

\medskip

\noindent{This law is valid up to the frequency $ k_d $ with}

\medskip

\noindent{\cite{L}}-VI-4-2 $ \qquad \, k_d \ \sim  \ ( \, e \, / \, \kappa^3 \, )^{1/4} \, . $

\medskip

\noindent{The small quantity $ \eps := k_d^{-1} $ is the Kolmogorov 
dissipative scale.} The relations \cite{L}}-VI-3-15 and \cite{L}}-VI-4-2 
imply that the rate of injection of kinetic energy $ e $ is linked to 
the number $ l $ according to $ e \sim \eps^{-1+\frac{3}{l}} $. We 
recover here that $ e \sim 1 $ when $ l =3 $ (see \cite{C}).

\medskip

\noindent{A starting point for the conventional theory of turbulence is
the notion that, on average, kinetic energy is transfered from low wave
numbers modes to high wave numbers modes. A recent paper \cite{FMRT} put 
forward the following idea: in the spectral region below that of injection 
of energy, an inverse (from high to low modes) transfer of energy takes 
place. At any rate, it is a central question to determine how the kinetic
energy is distributed.

\subsection{Phenomenological comparison.} The statistical 
approach is concerned with the spectral properties of solutions.
Below, we draw a parallel with the propagation of quasi-singularities
as it is described in the Theorem \ref{appBKW}. 

\medskip

\noindent{Let us examine how the square $ \cF(\u^\eps_\flat) (t,\xi)^2 $ 
of the Fourier transform of $ \u^\eps_\flat(t,x) $ is distributed.} To 
this end, consider the application
\vskip -5mm
$$ \left. \begin{array} {rcl}
\tilde E(t,\cdot) \, : \, \RR^+ & \longrightarrow & \RR^+ \\
k \ & \longmapsto & \tilde E(t,k) \, := \, \int_{\{ \xi 
\in \RR^d \, ; \, \vert \xi \vert = k \}} \ \vert \cF 
(\u^\eps_\flat) (t,\xi) \vert^2 \ \, d \sigma (\xi) \, . \quad
\end{array} \right. \ $$
\vskip -1mm
\noindent{The initial data $ \u^\eps_\flat (0,\cdot) $ has a {\it spectral gap}.} 
In another words, the graph of the function $ k \longmapsto \tilde
E(0,k) $ appears concentrated around the two characteristic wave 
numbers $ k \simeq 1 $ and $ k \simeq\eps^{-1} = k_d $.} In view of 
(\ref{phasenonn}), this situation does not persist. At the time 
$ t = \eps^{\frac{1}{l}} $, the concentration is around $ l $ 
characteristic wave numbers which are intermediate between the 
two preceding ones. This corresponds to a {\it discrete} cascade 
of energy.

\medskip

\noindent{Suppose now (\ref{speini}) and consider $ \dot \u^\eps_\flat $.} 
The life span of $ \dot \u^\eps_\flat (t,\cdot) $ is $ \eps^{\frac{1}{l}} 
\, T $.} There are various manners to get a family $ \{\dot \u^\eps_\flat(t,\cdot)
\}_{\eps \in \, ]0,1] } $ which is defined on some interval $ [0,\tilde T] $ with 
$ \tilde T > 0 $ independent on $ \eps $.} In particular, we can  

\medskip

\noindent{a) Select any $ \tilde T > 0 $ when $ T = + \infty $.} However 
nothing guarantees that the functions $ \dot \u^\eps_\flat $ are still 
approximate solutions on the interval $ [0,\tilde T] $. Indeed, since   
$ t $ is replaced by $ \eps^{- (1/l)} \, t $, the size of the error 
terms $ \dot \f^\eps_\flat $ depends on the increase of $ \f^\eps_\flat $  
with respect to $ t $. At this level, we are faced with secular growth 
problems \cite{La}.

\medskip

\noindent{b) Use a convergence process\footnote{When performing
the formal analysis, arbitrary values can be given to the 
parameters $ \eps \in \, ]0,1] $ and $ l \in \NN_* $. For instance 
$ \eps $ can be fixed whereas $ l $ goes to $ \infty $. Or 
$ l = - (\ln \eps) / (\ln 2) $ so that $  \eps^{\frac{1}{l}} \, 
T = \frac{1}{2} \, T > 0 $.} which needs the introduction of 
an {\it infinite} cascade of phase shifts. The intuition\footnote{
Even at a formal level, difficulties occur in order to justify 
the different convergences. Rigorous results in this direction
seem to be a difficult task.} is that the graph of $ \tilde E $ 
becomes continuous (no more gap). This corresponds to the 
impression of an {\it infinite} cascade of energy. This 
remark is consistent with engineering experiments and 
the observations reported in the statistical approach.

\bigskip

The turbulent phenomena which we study are very complex in their
realization. When $ t > 0 $, the description of $ \dot \u^\eps_\flat 
(t,\cdot) $ involves an infinite set of phases so that computations 
and representations are hard to implement. It gives the impression 
of a chaos. Nevertheless, our analysis reveals that these phenomena 
contain no mystery in their generation. On the contrary quantitative
and qualitative features can be predicted in the framework of non 
linear geometric optics. 

\smallskip

$\, $

\section{Euler equations in the variables $ (t,x,\theta) $.} 

As explained in the previous chapter, the demonstration of
the Theorem \ref{appBKW} is achieved with the representation
\begin{equation} \label{BKWdelbis}
\ \tilde \u^\eps_\flat (t,x) = \tilde u^\eps_\flat \bigl( 
t,x, \eps^{-1} \, \varphi^\eps_\flat (t,x) \bigr) \, , \quad 
\ \tilde \p^\eps_\flat (t,x) = \tilde p^\eps_\flat \bigl( 
t,x, \eps^{-1} \, \varphi^\eps_\flat (t,x) \bigr) \, . \ \
\end{equation}
Recall that the complete phase $ \varphi^\eps_\flat (t,x) $ is 
\begin{equation} \label{phasecomplete}
\varphi^\eps_\flat (t,x) \, = \, \varphi^\eps_g (t,x) + \eps \ 
\varphi^\eps_a (t,x) \, = \, \varphi_0 (t,x) + \hbox{$ \sum_{
k=1}^N $} \ \eps^{\frac{k}{l}} \ \varphi_k (t,x) \
\end{equation}
and that the profiles $ \tilde u^\eps_\flat ( t,x, \theta ) $ 
and $ \tilde p^\eps_\flat ( t,x, \theta ) $ have the form
\begin{equation} \label{feclate}
\left. \begin{array}{l}
\tilde u^\eps_\flat ( t,x, \theta ) \, = \, \u_0 (t,x) + 
\sum_{k=1}^N \, \eps^{\frac{k}{l}} \ \tilde U_k ( t,x, 
\theta) \, , \qquad \qquad \qquad \qquad \, \\
\tilde p^\eps_\flat ( t,x, \theta ) \, = \, \p_0 (t,x) + 
\sum_{k=1}^N \, \eps^{\frac{k}{l}} \ \tilde P_k ( t,x, 
\theta) \, . 
\end{array} \right. 
\end{equation}

\subsection{Preliminaries.}

$ \bullet $ {\bf Anisotropic viscosity.} Mark the 
abbreviated notations

\medskip

$ X^\eps_\flat (t,x) := \nabla \varphi^\eps_\flat (t,x) = 
\sum_{k=0}^N \, \eps^{\frac{k}{l}} \ X_k (t,x) \, , \qquad 
X_k (t,x) := \nabla \varphi_k (t,x) \, , $

\medskip

$ \mX^\eps_{\flat 1} (t,x) := \vert X^\eps_\flat (t,x) 
\vert^{-1} \ X^\eps_\flat (t,x) \, . $

\medskip

\noindent{Complete the unit vector $ \mX^\eps_{\flat 1} (t,x) $  
into some orthonormal basis of $ \RR^d $} 

\medskip

$ \mX^\eps_{\flat i} (t,x) \cdot \mX^\eps_{\flat j} (t,x) =
\delta_{ij} \, , \qquad \forall \, (i,j) \in \{1, \cdots , 
d \}^2 \, , $

\medskip

\noindent{so that all the vector fields $ \mX^\eps_{\flat i} $
are smooth functions on $ [0,T] \times \RR^d $.} The 
corresponding differential operators are denoted

\medskip

$ \mX^\eps_{\flat i} (\part ) := \mX^\eps_{\flat i} (t,x) 
\cdot \nabla \, , \qquad i \in \{1, \cdots , d \} \, . $

\medskip

\noindent{Their adjoints are} 

\medskip

$ \mX^\eps_{\flat i} (\part )^* := \mX^\eps_{\flat i} (t,x) 
\cdot \nabla + \Div \, ( \mX^\eps_{\flat i} ) (t,x) \, , 
\qquad i \in \{1, \cdots , d \} \, . $

\medskip

\noindent{Select $ \mq \in C^\infty_b ( [0,T] \times \RR^d ;
S^d_+ ) $ be such that}

\medskip

$ \exists \, c > 0 \, ; \qquad \mq (t,x) \geq c \, , \qquad 
\forall \, (t,x) \in [0,T] \times \RR^d \, . $

\medskip

\noindent{Let $ (m,n ) \in \NN^2 $. Consider the elliptic 
operator $ E^{\eps m}_{\flat n} (\part) $ defined according
to}

\medskip

$ E^{\eps m}_{\flat n} (\part) \, := \,  \bigl( \,
\eps^{\frac{m}{l}} \, \mX^\eps_{\flat 1} (\part )^* \, ,
\, \eps^{\frac{n}{l}} \, \mX^\eps_{\flat 2} (\part )^* \, ,
\, \cdots \, , \, \eps^{\frac{n}{l}} \, \mX^\eps_{\flat d} 
(\part )^* \, \bigr) $
$$ \qquad \qquad \qquad \qquad 
\left( \begin{array} {cccc}
\mq_{11} (t,x) & \mq_{12} (t,x) &\cdots & \mq_{1d} (t,x) \\
\mq_{21} (t,x) & \mq_{22} (t,x) & \cdots & \mq_{2d} (t,x) \\
\vdots & \vdots & \ & \vdots \\
\mq_{d1} (t,x) & \mq_{d2} (t,x) & \cdots & \mq_{dd} (t,x) 
\end{array} \right)
\left( 
\begin{array} {c}
\eps^{\frac{m}{l}} \, \mX^\eps_{\flat 1} (\part ) \\
\eps^{\frac{n}{l}} \, \mX^\eps_{\flat 2} (\part ) \\
\vdots \\
\eps^{\frac{n}{l}} \, \mX^\eps_{\flat d} (\part )
\end{array} \right) . $$
The introduction of the operator $ E^{\eps m}_{\flat n} (\part) $
in the right of $ (\cE) $ is compatible with the propagation
of oscillations only if $ m \geq l $ and $ n \geq 0 $. We
retain the limit case $ l=m $ and $ n = 0 $. The other 
situations are easier to deal with, at least when performing 
formal computations.

\bigskip

\noindent{$ \bullet $ {\bf Interpretation in $ (t,x,\theta) $.}} 
To deal with the variables $ (t,x,\theta) $, define

\medskip

$ \md_{j,\eps} \, := \, \eps \ \part_j \, + \, \part_j 
\varphi^\eps_\flat \times \part_\theta \, , \qquad j \in 
\{0, \cdots, d \} \, , $

\medskip

$ \md_\eps \, := \, ( \md_{1,\eps} , \cdots , \md_{d,\eps} ) 
\, , $

\medskip

$ \Grade^\eps_\flat \, := \, {}^t ( \md_{1,\eps} , \cdots ,
\md_{d,\eps} ) \, = \, \eps \ \nabla + X^\eps_\flat \times
\partial_\theta \, , $

\medskip

$ \Dive^\eps_\flat \, := \, ( \Grade^\eps_\flat)^\star \,
= \, \eps \ \Div + X^\eps_\flat \cdot \partial_\theta \, . $

\medskip

\noindent{The derivatives $ \mX^\eps_{\flat \dag} $ become}

\medskip

$ \eps \ \mX^\eps_{\flat 1} ( \md_\eps ) \, := \, \eps \
\mX^\eps_{\flat 1} (\part ) \, + \, \vert X^\eps_\flat 
(t,x) \vert \times \part_\theta \, , $

\medskip

$ \eps \ \mX^\eps_{\flat j} ( \md_\eps ) \, := \, \eps \
\mX^\eps_{\flat j} (\part ) \, , \qquad \forall \, j \in 
\{2, \cdots, d \} \, . $

\medskip

\noindent{The action of $ E^{\eps l}_{\flat 0} (\part) $ 
expressed in the variables $ (t,x,\theta) $ gives rise to
some negative differential operator of the order two, noted
$ E^{\eps l}_{\flat 0} ( \md_\eps ) $.} The coefficients
of the derivatives in $ E^{\eps l}_{\flat 0} ( \md_\eps ) $
are of size one, except in front of $ \mX^\eps_{\flat 1} 
(\part ) $. To avoid technicalities and to simplify the
notations, we substitute the Laplacian $ \nu \, \Delta $
for $ E^{\eps l}_{\flat 0} ( \md_\eps ) $. 

\medskip

\noindent{When $ \nu = 0 $, we recover Euler equations.} 
When $ \nu > 0 $, the action $ \nu \, \Delta $ can be 
viewed as the `trace' in $ (t,x,\theta) $ of the anisotropic 
viscosity $ E^{\eps l}_{\flat 0} (\part) $. Now, consider 
the Cauchy problem
\begin{equation} \label{cauchypr}
\ \left \lbrace \begin{array} {ll}
\md_{0,\eps} \, \tilde u^\eps_\flat + ( \tilde u^\eps_\flat
\cdot \Grade^\eps_\flat) \, \tilde u^\eps_\flat \! \! \! & 
+ \, \Grade^\eps_\flat \, \tilde p^\eps_\flat \\
\ & = \, \nu \ \eps \ \Delta \, \tilde u^\eps_\flat + \tilde 
f^\eps_\flat \, , \qquad \Dive^\eps_\flat \, \tilde 
u^\eps_\flat = \tilde g^\eps_\flat \, , \qquad \\ 
\tilde u^\eps_\flat (0,\cdot) \, = \, \tilde h^\eps_\flat 
(\cdot) \, , & \
\end{array} \right. 
\end{equation}
with given data 

\medskip

$ \tilde f^\eps_\flat \in \cW^\infty_T \, , \qquad \tilde 
g^\eps_\flat \in \cW^\infty_T \, , \qquad \tilde h^\eps_\flat 
\in H^\infty \, . $ 

\medskip

\noindent{Suppose that $ \nu = 0 $ and select some smooth 
solution $ ( \tilde u^\eps_\flat , \tilde p^\eps_\flat ) $ 
of (\ref{cauchypr}). The expressions $ \tilde \u^\eps_\flat $ 
and $ \tilde \p^\eps_\flat $ given by the formula (\ref{BKWdelbis}) 
are subjected to
\begin{equation} \label{cauchyprbf}
\left \lbrace \begin{array} {l}
\part_t \tilde \u^\eps_\flat + ( \tilde \u^\eps_\flat \cdot \nabla) 
\tilde \u^\eps_\flat + \nabla \tilde \p^\eps_\flat \, = \, \tilde 
\f^\eps_\flat \, , \qquad \Div \, \tilde \u^\eps_\flat = \tilde 
\g^\eps_\flat \, , \qquad \ \, \\
\tilde \u^\eps_\flat (0,\cdot) \, = \, \tilde \h^\eps_\flat (\cdot) \, , 
\end{array} \right. \qquad 
\end{equation}
where the functions $ \tilde \f^\eps_\flat (t,x) $, $ \tilde 
\g^\eps_\flat (t,x) $ and $ \tilde \h^\eps_\flat (t,x) $ are 
obtained by replacing the variable $ \theta $ by $ \varphi^\eps_\flat 
(t,x) $ in the expressions $ \eps^{-1} \ \tilde f^\eps_\flat 
( t,x, \theta ) $, $ \eps^{-1} \ \tilde g^\eps_\flat ( t,x, 
\theta ) $ and $ \tilde h^\eps_\flat (t,x, \theta) $. In other 
words, any solution of (\ref{cauchypr}) with $ \nu = 0 $ yields 
a solution of (\ref{cauchyprbf}). From now on, we proceed directly 
with the relaxed system (\ref{cauchypr}).

\smallskip

\subsection{The BKW analysis.}

Select a smooth solution $ \u_0(t,x) \in \cW^\infty_T $ of

\medskip

$ \part_t \u_0 + (\u_0 \cdot \nabla ) \, \u_0 + \nabla \p_0 \,
= \, \nu \ \Delta_x \u_0 \, , \qquad \Div \ \u_0 = 0 \, . $

\medskip

\noindent{Choose a phase $ \varphi_0 (t,x) \in C^1 ([0,T] 
\times \RR^d) $ with $ \nabla \varphi_0 (t,x) \in C^\infty_b 
([0,T] \times \RR^d) $. Suppose moreover that it satisfies 
the eiconal equation $ (ei) $ and the condition (\ref{nonstaou}).} 
The main step in the construction of approximate solutions is the 
following intermediate result.

\smallskip

\begin{prop} \label{appBKWinter} Select any $ \flat = (l,N) \in \NN^2 $ 
such that $ 0 < l < N $. Consider the following initial data

\medskip

$ \tilde U_{k0}^* (x,\theta) = \Pi_0 (0,x) \, \tilde U_{k0}^* 
(x,\theta) \in H^\infty \, , \qquad  1 \leq k \leq N \, , $

\smallskip

$ \langle \tilde U_{k0} \rangle (x) \in H^\infty \, , 
\qquad 1 \leq k \leq N \, , $

\smallskip

$ \varphi_{k0} (x) \in H^\infty \, , \qquad  \quad \!
1 \leq k \leq N \, . $

\medskip

\noindent{There are finite sequences $ \{ \tilde U_k \}_{1 \leq 
k \leq N} $ and $ \{ \tilde P_k \}_{1 \leq k \leq N} $ with}

\medskip

$ \tilde U_k (t,x,\theta) \in \cW^\infty_T \, , \qquad \tilde 
P_k (t,x,\theta) \in \cW^\infty_T \, , \qquad 1 \leq k \leq N \, , $

\medskip

\noindent{and a finite sequence $ \{ \varphi_k \}_{1 \leq k \leq N} $ 
with}

\medskip

$ \varphi_k (t,x) \in \cW^\infty_T \, , \qquad 1 \leq k \leq N \, , $

\medskip

\noindent{which are such that}

\medskip

$ \Pi_0 (0,x) \, \tilde U_k^* (0,x,\theta) =  \Pi_0 (0,x) \, 
\tilde U_{k0}^* (x,\theta) \, , \qquad 1 \leq k \leq N \, , $

\smallskip

$ \langle \tilde U_k \rangle (0,x) = \langle \tilde U_{k0} 
\rangle (x) \, , \qquad \! 1 \leq k \leq N \, , $

\smallskip

$ \varphi_k (0,x) = \varphi_{k0} (x) \, , \qquad \quad
\ 1 \leq k \leq N \, . $ 

\medskip

\noindent{Define $ \varphi^\eps_\flat $ as in (\ref{phasecomplete}).}
All the preceding expressions are adjusted so that the functions 
$ \tilde u^\eps_\flat $ and $ \tilde p^\eps_\flat $ associated with 
the expansions in (\ref{feclate}) are approximate solutions 
on the interval $ [0,T] $. More precisely, they satisfy  
(\ref{cauchypr}) with 
\begin{equation} \label{conclufi}
\tilde h^\eps_\flat (x, \theta) \, = \, \u_0 (0,x) + \hbox{$ \sum_{
k=1}^N $} \, \eps^{\frac{k}{l}} \ \tilde U_{k0} (x,\theta) \qquad 
\qquad \qquad
\end{equation}
and we have
\begin{equation} \label{conclufiaussi}
\{ \tilde f^\eps_\flat \}_\eps = \bigcirc( \eps^{\frac{N+1}{l}}) 
\, , \qquad \{ \tilde g^\eps_\flat \}_\eps = \bigcirc( \eps^{
\frac{N+1}{l}}) \, . \qquad \qquad \qquad \quad
\end{equation}

\end{prop}

\medskip

\noindent{$ \bullet $ {\bf Proof of the Proposition \ref{appBKWinter}.}}
For convenience, we will drop in this paragraph the tilde '$\, \tilde \ \, $' 
on the profiles $ u^\eps_\flat $, $ p^\eps_\flat $, $ U_k $ and $ P_k $.
This modification concerns only this demonstration. We hope that it 
will not induce confusions: we still work here with the complete phase 
$ \varphi^\eps_\flat $.
 
\medskip

\noindent{Because of (\ref{nondegeat}) we can define the application 
$ \Pi^\eps_\flat (t,x) $ which is the orthogonal projector on the 
hyperplane $ \nabla \varphi^\eps_\flat (t,x)^\perp \subset \RR^d $.} 
We adopt the convention

\medskip

$ \Pi^\eps_\flat (t,x) = \sum_{k=0}^\infty \, \eps^{\frac{k}{l}} \ 
\Pi_k (t,x) \, , \qquad \Pi_k \in \cW^\infty_T \, , \qquad \eps 
\in \, ]0, \eps_0] \, . $ 

\medskip

\noindent{The access to $ \Pi_k $ needs only the knowledge of the $ X_j $ 
for $ j \leq k $.} Introduce

\medskip

$ v^\eps_\flat := X^\eps_\flat \cdot u^\eps_\flat = \sum_{k=0}^\infty \,
\eps^{\frac{k}{l}} \ V_k \, , \qquad \! V_k = X_k \cdot \u_0 + 
\sum_{j=0}^{k-1} \, X_j \cdot U_{k-j} \, , $

\medskip

$ w^\eps_\flat := \Pi^\eps_\flat \, u^\eps_\flat = \sum_{k=0}^\infty \, 
\eps^{\frac{k}{l}} \ W_k \, , \qquad W_k = \Pi_k \, \u_0 + \sum_{j=0}^{k-1} 
\, \Pi_j \, U_{k-j} \, . $

\medskip

\noindent{By construction}

\medskip

$ u^\eps_\flat \, = \, v^\eps_\flat \ \vert X^\eps_\flat \vert^{-2} 
\ X^\eps_\flat + w^\eps_\flat \, , \qquad U_k \, = \, V_k \ \vert X_0 
\vert^{-2} \  X_0 \, + W_k + h $

\medskip

\noindent{where $ h $ depends only on the $ X_j $ for $ j \leq k $ 
and on the $ U_j $ for $ j \leq k-1 $.} 

\medskip

\noindent{The conditions prescribed in the Proposition \ref{appBKWinter} 
on the initial data $ U_{k0}^* $ allow to fix the functions $ \nabla 
\varphi_0(0,x) \cdot U_k^* (0,x,\theta) $ as we want.} Since

\medskip

$ V_k^* \, = \, \nabla \varphi_0 \cdot U_k^* \, + \, \sum_{j=1}^{k-1} \, 
X_j \cdot U_{k-j}^* \, , $ 

\medskip

\noindent{the same is true (by induction) for the components $ V_k^* 
(0,x,\theta) $.} To begin with, we impose the polarization conditions
\begin{equation} \label{polafirst}
P_k^* \, \equiv \, V_k^* \, \equiv \, 0 \, , \qquad \forall \, k \in 
\{ 1,\cdots,l \} \qquad \qquad \qquad \qquad \
\end{equation}
and we adjust {\it a priori} the geometrical phase $ \varphi^\eps_g $  
so that
\begin{equation} \label{shiftfirst}
\part_t \varphi_k + \bar V_k = 0 \, , \qquad \forall \, k \in \{ 1,
\cdots, l-1 \} \qquad \qquad \qquad \qquad \quad
\end{equation}
which implies that

\medskip

$ \part_t \varphi^\eps_g + (\bar u^\eps_\flat \cdot \nabla)
\varphi^\eps_g \, = \, \sum_{k=l}^\infty \, \eps^{\frac{k}{l}} 
\ \bar V_k \, = \, \bigcirc (\eps) \, . $

\medskip

\noindent{It amounts to the same thing to look at the equations in
(\ref{cauchypr}) or at the following singular system (we drop here 
the indices $ \eps $ and $ \flat $ at the level of $ u^\eps_\flat $, 
$ v^\eps_\flat $, $ w^\eps_\flat $, $ p^\eps_\flat $ , $ \Pi^\eps_\flat $
and $ \varphi^\eps_\flat $)}
\begin{equation} \label{singu}
\left \{ \begin{array}{ll}
\part_t u + ( u \cdot \nabla ) \, u + \nabla p & \! \! \! + \, 
\eps^{-1} \ (\part_t \varphi + v) \ \partial_\theta u \quad \\
& \! \! \! + \, \eps^{-1} \ \part_\theta p \ \nabla \varphi \, 
= \, \nu \ \Delta \, u \, , \\
\Div \ u + \eps^{-1} \ \part_\theta v = 0 \, . & \
\end{array} \right. \qquad \qquad \
\end{equation}
The functions $ v $ is subjected to
\begin{equation} \label{v}
\left. \begin{array}{rl}
\part_t v \! \! \! \! & + \, (u \cdot \nabla) \, v + X \cdot \nabla p 
+ \eps^{-1} \ (\part_t \varphi + v) \ \partial_\theta v  \\
\ & + \, \eps^{-1} \ \part_\theta p \ \parallel X \parallel^2 
- \bigl( \part_t X + (u \cdot \nabla) \, X \bigr) \cdot u  = 
\nu \  X \cdot \Delta \, u \, .
\end{array} \right.
\end{equation}
The functions $ w $ satisfies
\begin{equation} \label{w}
\left. \begin{array}{rl}
\part_t w \! \! \! \! & + \, (u \cdot \nabla) \, w + \Pi \, \nabla 
p + \eps^{-1} \ (\part_t \varphi + v) \ \partial_\theta w \qquad 
\qquad \quad \ \ \ \\
\ & - \, \bigl( \part_t \Pi + (u \cdot \nabla) \, \Pi \bigr) \, u 
= \nu \ \Pi \, \Delta \, u \, .
\end{array} \right.
\end{equation}
Substitute the expressions $ u^\eps_\flat $ and $ p^\eps_\flat $ 
given by (\ref{feclate}) into (\ref{singu}). Then arrange the terms 
according to the different powers of $ \eps $ which are in factor. 
The contributions coming from the orders $ \eps^{\frac{1}{l}-1} $, 
$ \cdots $, $ \eps^{-\frac{1}{l}} $ and $ \eps^0 $ are eliminated 
through (\ref{polafirst}), (\ref{shiftfirst}) and the constraints 
imposed on $ (\u_0,\p_0) $. 

\medskip

\noindent{Now, look at the terms in front of $\eps^{\frac{j}{l}} $ 
with $ j \in \NN_* $.} It remains
\begin{equation} \label{eqgj}
\left \{ \begin{array}{l}
\part_t U_j + \sum_{k=0}^j \, ( U_k \cdot \nabla) \, U_{j-k} + 
\nabla P_j + \sum_{k=0}^j \, \part_\theta P_{l+k} \ \nabla 
\varphi_{j-k} \\
\qquad \, + \, \sum_{k=0}^{j-1} \, (\part_t \varphi_{l+k} + V_{l+k}) \ 
\part_\theta U_{j-k} \, = \, \nu \ \Delta \, U_j \, , \\
\Div \ U_j + \part_\theta V_{j+l} = 0 \, .
\end{array} \right. \quad
\end{equation}
Proceed in a similar manner with (\ref{v}). Just arrange the terms
which come from $ X \cdot \Delta \, u $ and which do not involve 
$ U_j $ in a source term $ \cH^j_V $
\begin{equation} \label{eqgjeg}
\left. \begin{array}{rl}
\, \part_t V_j \! \! \! \! & + \, \sum_{k=0}^j \, (U_k \cdot \nabla) \, 
V_{j-k} + \sum_{k=0}^{j-1} \, (\part_t \varphi_{l+k} + V_{l+k}) \
\part_\theta V_{j-k} \\
\ & + \, \sum_{k=0}^j \, X_k \cdot \nabla P_{j-k} \, - \, \sum_{k=0}^j \, 
\part_t X_k \cdot U_{j-k} \\
\ & - \, \sum_{k=0}^j \, \bigl( \sum_{l=0}^k \, (U_{k-l} \cdot \nabla) \, 
X_l \bigr) \cdot U_{j-k} \\
\ & + \, \sum_{k=1}^{j-1} \, \bigl( \sum_{l=0}^k \, X_l \cdot X_{k-l} 
\bigr) \ \part_\theta P_{j+l-k} + \vert X_0 \vert^2 \
\part_\theta P_{j+l} \\
\ & = \, \nu \ X_0 \cdot \Delta \, U_j + \cH^j_V (t,x,\theta,U_1,X_1, 
\cdots, U_{j-1},X_{j-1}, X_j) \, . \ \
\end{array} \right.
\end{equation}
The same operation with (\ref{w}) yields 
\begin{equation} \label{eqgjel}
\, \left. \begin{array}{rl}
\part_t W_j \! \! \! \! & + \, \sum_{k=0}^j \, (U_k \cdot \nabla) \, 
W_{j-k} \, + \, \sum_{k=0}^{j-1} \, (\part_t \varphi_{l+k} + V_{l+k}) \ 
\part_\theta W_{j-k} \quad \\
\ & + \, \sum_{k=0}^j \, \Pi_k \, \nabla P_{j-k} \, - \, \sum_{k=0}^j \, 
\part_t \Pi_k \, U_{j-k} \\
\ & - \, \sum_{k=0}^j \, \bigl( \sum_{l=0}^k \, (U_{k-l} \cdot \nabla)
\, \Pi_l \bigr) \,U_{j-k} \\ 
\ & = \, \nu \ \Pi_0 \ \Delta \, U_j + \cH^j_W (t,x,\theta,U_1,X_1, 
\cdots, U_{j-1},X_{j-1}, X_j) \, .
\end{array} \right.
\end{equation}
Then extract the mean value of (\ref{eqgj}) 
\begin{equation} \label{moyj}
\ \left \{ \begin{array}{l}
\part_t \bar U_j + (\u_0 \cdot \nabla) \, \bar U_j + (\bar U_j \cdot 
\nabla) \, \u_0 + \nabla \bar P_j \qquad \qquad \qquad \qquad \qquad \\
\qquad \, + \, \sum_{k=1}^{j-1} \, \langle ( U_k \cdot \nabla) \, 
U_{j-k} \rangle \\
\qquad \, + \sum_{k=1}^{j-1} \, \langle V^*_{l+k} \ \part_\theta 
U_{j-k} \rangle \, = \, \nu \ \Delta_x \, \bar U_j \, , \\
\Div \ \bar U_j = 0 \, .
\end{array} \right.
\end{equation}
Observe also that
\begin{equation} \label{divj}
V^*_{j+l} \, = \, - \, \Div \ \partial_\theta^{-1} U^*_j \, , \qquad 
\forall \, j \in \NN_* \, . \qquad \qquad \qquad \qquad \qquad
\end{equation}
Consider the inductive reasoning based on

\medskip

\noindent{\bf Hypothesis} $ (H_j) \, $:

\medskip

{\em
\noindent{i) The expressions $ U_1 $, $ \cdots $, $ U_j $ and  $ P_1 $, 
$ \cdots $, $ P_j $ are known.}

\medskip

\noindent{ii) The phases $ \varphi_1 $, $\cdots $, $ \varphi_j $ are 
identified. The same is true for the vectors $ X_1 $, $\cdots $, $ X_j $ 
and the projectors $ \Pi_1 $, $\cdots $, $ \Pi_j $. Moreover, the 
following relations are satisfied}
\begin{equation} \label{phik}
\part_t \varphi_{j+k} + \bar V_{j+k} = 0 \, , \qquad \forall \, k \in 
\{ 1, \cdots, l-1 \} \, . 
\qquad \end{equation}
iii) The correctors $ V^*_{j+1} $, $ \cdots $, $ V^*_{j+l} $ and 
$ P^*_{j+1} $, $ \cdots $, $ P^*_{j+l} $ are identified and
\begin{equation} \label{vk+n}
V^*_{j+k} = \, - \, \Div \ \part_\theta^{-1} U^*_{j+k-l} \, , \qquad 
\forall \, k \in \{1, \cdots, l \}  \, . \qquad \quad
\end{equation}
}

\noindent{- \underline{Verification} \underline{of} $ (H_1) $}. The mean
value $ \bar U_1 $ is obtained by solving 

\medskip

$ \part_t \bar U_1 + (\u_0 \cdot \nabla) \, \bar U_1 + (\bar U_1 \cdot 
\nabla) \, \u_0 + \nabla \bar P_1 = \nu \ \Delta_x \, \bar U_1 \, , 
\qquad \Div \ \bar U_1 = 0 \, . $

\medskip

\noindent{Using (\ref{shiftfirst}), it allows to determine $ \varphi_1 $.}
Now look at the oscillating part of (\ref{eqgjel}) with the indice $ j = 1 $.
The constraint on $ W_1^* \equiv U_1^* $ writes

\medskip

$ \part_t W_1^* + ( \u_0 \cdot \nabla) \, W_1^* + ( \part_t \varphi_l 
+ \bar V_l) \ \part_\theta W_1^* = M \, W_1^* + \nu \ \Pi_0 \ \Delta
\, W^*_1 $

\medskip

\noindent{where $ M $ is the linear application}

\medskip

$ M \, U \, := \, (\part_t \Pi_0) \, U + \bigl( (\u_0 \cdot \nabla) 
\Pi_0 \bigr) \, U - \Pi_0 \, (U \cdot \nabla ) \u_0 \, . $

\medskip

\noindent{We impose (\ref{phik}) for $ j=1 $.} In view of 
(\ref{shiftfirst}), it reduces to

\medskip

$ \part_t \varphi_l + \bar V_l = 0 \, . $

\medskip

\noindent{The link between $ W_1^* $ and $ \bar V_l $ is 
removed.} It remains the linear equation

\medskip

$ \part_t W_1^* + ( \u_0 \cdot \nabla) \, W_1^* = M \, W_1^* + \nu \ 
\Pi_0 \ \Delta \, W^*_1 \, . $

\medskip

\noindent{When $ \nu = 0 $, the profile $ W^*_1 $ is obtained 
by integrating along the characteristic curves $ \Gamma (\cdot,x) $.
This justifies the remark 2.2.4 for $ k = 1 $. Observe also that the 
polarization condition $ W^*_1 = \Pi_0 \, W^*_1 $ is conserved since 
the equation given for $ W_1^* $ is equivalent to}
\begin{equation} \label{eqW} 
\ \left \{ \begin{array}{l}
\Pi_0 \, \bigl \lbrack \part_t W_1^* + ( \u_0 \cdot 
\nabla) \, W_1^* + (W_1^* \cdot \nabla) \u_0  \bigr \rbrack \, 
= \, \nu \ \Pi_0 \ \Delta \, W^*_1 \, , \qquad  \quad \\
W_1^* \, = \, \Pi_0 \, W^*_1 \, .   
\end{array} \right. 
\end{equation}
Introduce the linear form 

\medskip

$ \ell \, U \, := \ \vert X_0 \vert^{-2} \ \bigl \lbrack \, 
\part_t X_0 \cdot U + \bigl( (\u_0 \cdot \nabla) X_0 \bigr) 
\cdot U - X_0 \cdot \bigl( (U \cdot \nabla ) \u_0 \bigr) \, 
\bigr \rbrack \, . $

\medskip

\medskip

\noindent{The constraint $ V_1^* \equiv 0 $ is equivalent to}

\medskip

$ P^*_{l+1} \, = \, \ell \ \part_\theta^{-1} W_1^* \, + \, 
\nu \ \vert X_0 \vert^{-2} \ X_0 \cdot \part_\theta^{-1} \, 
\Delta \, W^*_1 \, , $ 

\medskip

\noindent{We have also}

\medskip

$ V^*_{l+1} = \, - \, \Div \ \part_\theta^{-1} W_1^* \, . $

\medskip

\noindent{At this stage, we know who is $ U_1 \equiv \bar U_1 
+ W_1^* $ and $ P_1 \equiv 0 $.} Moreover, we have the relations 
(\ref{phik}) and (\ref{vk+n}). Thus, the hypothesis $ (H_1) $ 
is verified.

\bigskip

\noindent{- \underline{The} \underline{induction}.} Suppose 
that the conditions given in $ (H_j )$ are satisfied. The 
question is to obtain $ (H_{j+1} )$. Consider first (\ref{moyj}) 
with the indice $ j+1 $. The relation (\ref{vk+n}) induces 
simplifications. It remains
\begin{equation} \label{moyj+1}
\left \{ \begin{array}{l}
\part_t \bar U_{j+1} + (\u_0 \cdot \nabla) \, \bar U_{j+1} 
+ (\bar U_{j+1} \cdot \nabla) \, \u_0 \qquad \qquad \\
\qquad + \, \sum_{k=1}^{j} \, ( \bar U_k \cdot \nabla) \, 
\bar U_{j+1-k} + \, \sum_{k=1}^j \, \Div \ \langle 
U^*_k \otimes U^*_{j+1-k} \rangle \\
\qquad + \, \nabla \bar P_{j+1} \, = \, \nu \ \Delta_x 
\, \bar U_{j+1} \, , \qquad \Div \ \bar U_{j+1} = 0 \, .
\end{array} \right.
\end{equation}
This system gives access to $ \bar U_{j+1} $ and $ \bar P_{j+1} $. 
For $ j=1 $, it yields
\begin{equation} \label{moy2}
\ \left \{ \begin{array}{l}
\part_t \bar U_2 + (\u_0 \cdot \nabla) \, \bar U_2 + (\bar U_2 \cdot 
\nabla) \, \u_0 + \nabla \bar P_2 \qquad \qquad \qquad \quad \ \ \, \\
\quad \ + \,( \bar U_1 \cdot \nabla) \, \bar U_1 + \Div \ \langle 
U^*_1 \otimes U^*_1 \rangle = \, \nu \ \Delta \bar U_2 \, , \quad
\ \Div \ \bar U_2 = 0 \, .
\end{array} \right.
\end{equation}
Because of (\ref{nonnontri}), the source term $ \langle U^*_1 \otimes
U^*_1 \rangle $ is sure to be non trivial. We recover here that in
general $ \bar U_2 \not \equiv 0 $ even if $ \bar U_1 (0,\cdot) 
\equiv \bar U_2  (0,\cdot) \equiv 0 $. The term $ \bar U_2 $ 
excites $ \varphi_2 $ through (\ref{trans2}). Generically, we 
have $ \varphi_2 \not \equiv 0 $ even if

\medskip

$ \bar U_1 (0,\cdot) \, \equiv \, \bar U_2  (0,\cdot) \, \equiv \, 0 \, , 
\qquad \varphi_1 (0,\cdot) \, \equiv \, \varphi_2  (0,\cdot) \, \equiv 
\, 0 \, . $

\medskip

\noindent{Observe however that exceptions can happen (see the 
remark 3.5.4).} The information (\ref{phik}) for $ k=1 $ means 
that

\medskip

$ \part_t \varphi_{j+1} + (\u_0 \cdot \nabla) \, \varphi_{j+1} 
+ X_0 \cdot \bar U_{j+1} + \sum_{l=1}^j \, X_l \cdot \bar U_{j
+1-l} \, = \, 0 \, . $

\medskip

\noindent{Deduce $ \varphi_{j+1} $ from this equation, and 
therefore $ X_{j+1} $ and $ \Pi_{j+1} $.} Complete with the 
triangulation condition
\begin{equation} \label{phibi}
\part_t \varphi_{j+l} + \bar V_{j+l} = 0 \, . \qquad \qquad \qquad
\qquad \qquad \qquad \qquad \qquad \quad \quad
\end{equation}
Then extract the oscillating part of (\ref{eqgjel}) written 
with $ j+1 $. Use $ (H_j) $ and (\ref{phibi}) in order to 
simplify the resulting equation. It yields
\begin{equation} \label{oscij+1}
\ \part_t W_{j+1}^* + ( \u_0 \cdot \nabla) \, W_{j+1}^* = M \, 
W_{j+1}^* + \nu \ \Pi_0 \ \Delta \, W^*_{j+1} + f \qquad \qquad \
\end{equation}
where $ f $ is known. We get $ W^*_{j+1} $ by solving (\ref{oscij+1}).
Therefore we have $ U^*_{j+1} $ and we can deduce $ V^*_{j+l+1} = \, 
- \, \Div \ \part_\theta^{-1} U^*_{j+1} $.

\medskip

\noindent{Now look at the constraint (\ref{eqgjeg}) for the 
indice $ j+1 $.} Extract the oscillating part. It allows to 
recover $ P^*_{j+l+1} $. Thus we have $ (H_{j+1}) $. 

\medskip

\noindent{Apply the induction up to $ j = N -l $.} It yields 
$ U_1 $, $ \cdots $, $ U_N $. Construct oscillations 
$ \tilde u^\eps_\flat $ and $ \tilde p^\eps_\flat $ by way of 
(\ref{feclate}). It furnishes source terms $ \tilde f^\eps_\flat $ 
and $ \tilde g^\eps_\flat $ through (\ref{cauchypr}). By 
construction, we recover (\ref{conclufiaussi}). \hfill $ \Box $

\subsection{Divergence free approximate solutions in $ (t,x, \theta) $.} 

In this subsection 4.3, we impose on $ \varphi_0 $ a constraint 
which is more restrictive than (\ref{nonstaou}). We suppose that 
we can find a direction $ \zeta \in \RR^d \setminus \{ 0 \} $ such
that
\begin{equation} \label{nonstaoubis} 
\exists \, c > 0 \, ; \qquad \nabla \varphi_0 (t,x) \cdot \zeta \, 
\geq \, c \, , \qquad \forall \, (t,x) \in [0,T] 
\times \RR^d \, . \quad \
\end{equation}

\begin{prop} \label{complement} The assumptions are as in
the Proposition \ref{appBKWinter}. The profiles $  \tilde 
U_k $, $  \tilde P_k $, and the phases $  \varphi_k $ are
defined in the same way. Then, there are correctors

\medskip

$ \tilde {c u}^\eps_\flat (t,x,\theta) \in \cW^\infty_T \, , 
\qquad \{ \tilde{c u}^\eps_\flat \}_\eps \, = \, \bigcirc(
\eps^{\frac{N}{l}}) $ 

\medskip

\noindent{such that the functions $ \tilde u^\eps_\flat $ 
and $ \tilde p^\eps_\flat $ defined according to} 
$$ \, \left. \begin{array}{l}
\tilde u^\eps_\flat (t,x) := \u_0 (t,x) + \sum_{k=1}^N \, 
\eps^{\frac{k}{l}} \ \tilde U_k \bigl( t,x, \eps^{-1} \, 
\varphi^\eps_g (t,x) \bigr) + \tilde {c u}^\eps_\flat (t,x)  
\qquad \ \ \\
\tilde p^\eps_\flat (t,x) := \p_0 (t,x) + \sum_{k=1}^N \, 
\eps^{\frac{k}{l}} \ \tilde P_k \bigl( t,x, \eps^{-1} \, 
\varphi^\eps_g (t,x) \bigr) 
\end{array} \right. $$
satisfy the Cauchy problem
\begin{equation} \label{eclacpa}
\ \left \lbrace \begin{array} {l}
\md_{0,\eps} \, \tilde u^\eps_\flat + ( \tilde u^\eps_\flat
\cdot \Grade^\eps_\flat) \, \tilde u^\eps_\flat + 
\Grade^\eps_\flat \, \tilde p^\eps_\flat \\
\qquad \qquad \qquad \qquad \qquad \! = \, \nu \ \eps \ 
\Delta \, \tilde u^\eps_\flat + \tilde f^\eps_\flat \, , 
\qquad \Dive^\eps_\flat \, \tilde u^\eps_\flat = 0 \qquad \\ 
\tilde u^\eps_\flat (0,x,\theta) \, = \,  \u_0 (0,x) + 
\sum_{k=1}^N \, \eps^{\frac{k}{l}} \ \tilde U_{k0} (x,\theta)
\end{array} \right.
\end{equation}
and we still have $ \{ \tilde f^\eps_\flat \}_\eps \, = \, 
\bigcirc(\eps^{\frac{N+1}{l}}) $.

\end{prop}

\noindent{We need some material before proving the Proposition
\ref{complement}.}

\medskip

\noindent{$ \bullet $ {\bf The divergence free relation in 
the variables $ (t,x,\theta) $.}} We can select some special 
right inverse of the application $ \Dive^\eps_\flat \, : \, 
H^{\infty *}_T \longrightarrow H^{\infty *}_T $.} 

\begin{lem} \label{invpart*} There is a linear operator
$ \mr \mi \Dive^\eps_\flat \, : \, \text{\rm Im} \, ( 
\Dive^\eps_\flat) \longrightarrow H^{\infty *}_T $ with
\begin{equation} \label{exainv'}
\Dive^\eps_\flat \circ \mr \mi \Dive^\eps_\flat \ g \, = \, 
g \, , \qquad \forall \, g \in \text{\rm Im} \, (\Dive^\eps_\flat) 
\, . \qquad \qquad \qquad \qquad 
\end{equation}
For all $ m \in \NN $, there is a constant $ C_m > 0 $ such that
\begin{equation} \label{estiinv'}
\parallel \mr \mi \Dive^\eps_\flat \ g \parallel_{H^m} \, \leq \, 
C_m \ \parallel g \parallel_{H^{m+1+ \frac{d}{2}}} \, , 
\qquad \forall \, g \in \text{\rm Im} \, (\Dive^\eps_\flat) \, . 
\quad \ 
\end{equation}

\end{lem}

\medskip

\noindent{\em \underline{Proo}f \underline{o}f \underline{the} 
\underline{Lemma} 4\underline{.1}.} Let $ n \in \NN_* $. Note

\medskip

$ t_j := j \, T / n \, , \qquad x_j = k /n \, , \qquad 1 \leq 
j \leq n-1 \, , \qquad k \in \ZZ^d \, . $

\medskip

\noindent{Consider a related partition of unity}

\medskip

$ \chi_{(j,k)} \in C^\infty ([0,T] \times \RR^d) \, , \qquad 
\qquad (j,k) \in \{ 1, \cdots ,n-1 \} \times \ZZ^d \, , $

\smallskip

$ \sum_{j=1}^{n-1} \, \sum_{k \in \ZZ^d} \, \chi_{(j,k)} (t,x) \, 
= \, 1 \, , \qquad \forall \, (t,x) \in [0,T] \times \RR^d \, , $

\smallskip

$ \bigl \{ \, (t,x) \, ; \ \chi_{(j,k)}(t,x) \not = 0 \, \bigr \} \, 
\subset \, [t_j - \frac{2}{n} , t_j + \frac{2}{n} ] \times B (x_j,
\frac{2}{n} ] \, , $

\smallskip

$ \bigl \{ \, (t,x) \, ; \ \chi_{(j,k)}(t,x) = 1 \, \bigr \} \, 
\supset \, [t_j - \frac{1}{n} , t_j + \frac{1}{n} ] \times B 
(x_j,\frac{1}{n} ] \, . $

\medskip

\noindent{By hypothesis, there is a function $ v \in H^{\infty *}_T $ 
such that $ g = \Dive^\eps_\flat \, v $.} Introduce

\medskip

$ v_{(j,k)} \, := \, \chi_{(j,k)} \ v \in H^{\infty *}_T \, , 
\qquad g_{(j,k)} := \Dive^\eps_\flat \, v_{(j,k)} \, . $

\medskip

\noindent{It suffices to exhibit $ \mr \mi \Dive^\eps_\flat \, 
g_{(j,k)} $ and to show (\ref{estiinv'}) with a constant $ C_m $ 
which is uniform in $ {(j,k)} $.} The problem of finding $ \mr 
\mi \Dive^\eps_\flat \, g_{(j,k)} $ can be reduced to a model 
situation. This can be achieved by using a change of variables 
in $ (t,x) $, based on (\ref{nonstaoubis}). From now on, the 
time $ t $ is viewed as a parameter, the space variable is 
$ x = (x_1,\hat x) \in \RR \times \RR^{d-1} $, and we work with 

\medskip

$ g = g^* = \Dive^\eps_\flat \, v = ( \eps \, \part_1 + \part_\theta) 
v_1 + \part_2 v_2 + \cdots + \part_d v_d \, , $

\medskip

$ \bigl \{ \, x \, ; \ g(x,\theta) \not = 0 \, \bigr \} \, \subset \,
\bigl \{ \, x \, ; \ v(x,\theta) \not = 0 \, \bigr \} \, \subset \,
B(0,\frac{1}{2}] \, . $

\medskip

\noindent{Let $ \psi \in C^\infty (\RR^{d-1};\RR^+) $ be such that
$ \int_{\RR^{d-1}} \, \psi(\hat x) \ d \hat x = 1 $ and

\medskip

$ \bigl \{ \, \hat x \, ; \ \psi (\hat x) \not = 0 \, \bigr \} \, 
\subset \, B (0,1] \, , \qquad \bigl \{ \, \hat x \, ; \ \psi 
(\hat x) = 1 \, \bigr \} \, \supset \, B (0,\frac{1}{2}] \, . $

\medskip

\noindent{Decompose $ g $ according to}

$ g \, = \, (g - \breve g) \ \psi + \breve g \ \psi \, , \qquad 
\breve g (x) := \int_{\RR^{d-1}} \, g(x_1,\hat x) \ d \hat x \, 
= \, ( \eps \, \part_1 + \part_\theta) \breve v_1 \, . $

\medskip

\noindent{Seek a special solution $ u $ having the form}

\medskip

$ u \, = \, \mr \mi \Dive^\eps_\flat \, g \, = \, {}^t \bigl( \,
a \, , \, \text{\rm ridiv} \, [ (g - \breve g) \ \psi ] \, \bigr) \, , 
\qquad a \in H^{\infty *}_T $

\medskip

\noindent{where 'ridiv' is the operator of Lemma \ref{invpartbar}
applied in the dimension $ d-1 $. It remains to control the
scalar function $ a $ which satisfies the constraint}

\medskip

$ \eps \ \part_1 a + \part_\theta a \, = \, h \, := \, \breve g \ 
\psi \, = \, (\eps \ \part_1 + \part_\theta) (\breve v_1 \, 
\psi ) \, . $

\medskip

\noindent{Take the explicit solution}

\medskip

$ a (x_1, \hat x, \theta) \, = \, \int_{- \infty}^\theta \, h
\bigl( x_1 + \eps \, (s-\theta) , \hat x , s \bigr) \ ds $

\smallskip

$ \qquad \qquad \ \ = \, \eps^{-1} \ \int_{- \infty}^0 \, 
h \bigl(x_1 + r , \hat x , \theta + \eps^{-1} \, r \bigr) 
\ dr \, . $

\medskip

\noindent{By construction}

\medskip

$ a (x,\theta + 1) = a (x,\theta) \, , \qquad \int_\TT \, a 
(x,\theta) \ d \theta \, = \, 0 \, , \qquad \forall \, (x, 
\theta) \in \RR^d \times \TT \, . $

\medskip

\noindent{For $ \vert x_1 \vert + \vert \hat x \vert \geq 2 $, 
we find}

\medskip

$ a (x,\theta) \, = \, \int_{- \infty}^\theta \, \frac
{d}{ds} \, \bigl[ (\breve v_1 \, \psi) \bigl( x_1 + 
\eps \, (s-\theta) , \hat x , s \bigr) \bigr \rbrack \ 
ds \, = \, (\breve v_1 \, \psi ) (x,\theta) \, = \, 0 \, . $

\medskip

\noindent{It implies that}

\medskip

$ \bigl \{ \, (x,\theta) \, ; \ a (x,\theta) \not = 0 \, 
\bigr \} \, \subset \, B(0;2] \, . $

\medskip

\noindent{Note $ \mh := \part^{-1}_\theta h \in H^{\infty *} $.} 
Obviously

\medskip

$ \parallel \mh \parallel_{H^m} \, \leq \, C_m \ \parallel h 
\parallel_{H^m} \, , \qquad \forall \, m \in \NN \, , $

\medskip

$ \bigl \{ \, (x,\theta) \, ; \ \mh (x,\theta) \not = 0 \, 
\bigr \} \, \subset \, B(0;1] \, , $

\medskip

\noindent{and we have the identity}

\medskip

$ a (x_1,\hat x,\theta) \, = \, \mh (x_1,\hat x,\theta) \, 
- \, \int_{ - x_1 - 1}^{- x_1 + 1} \ \part_1 \mh \bigl(
x_1 + r , \hat x , \theta + \eps^{-1} \, r \bigr) \ dr \, . $

\medskip

\noindent{The term on the right is supported in $ B(0,2] $.}
Use Fubini and Cauchy-Schwarz inequality to control the 
integration of $ \part_1 \mh $. It yields (\ref{estiinv'}).} 
\hfill $ \Diamond $

\bigskip

\noindent{$ \bullet $ {\bf The Leray projector interpreted in 
the variables $ (t,x,\theta) $.}} Introduce the closed subspace

\medskip

$ \text{F}^\eps_\flat \, := \,  \bigl \lbrace \, u^* \in L^{2*}_T \, ; 
\ \Dive^\eps_\flat \, u^* = 0 \, \bigr \rbrace \, \subset \, 
L^{2*}_T \, . $

\medskip

\noindent{Note $ \mP^\eps_\flat $ the orthogonal projector from 
$ L^{2*}_T $ onto F$ {}^\eps_\flat $.} This is a self-adjoint 
operator such that

\medskip

$ \ker \, \Dive^\eps_\flat = \text{Im} \, \mP^\eps_\flat \, , \qquad 
\text{Im} \, \Grade^\eps_\flat = ( \ker \, \Dive^\eps_\flat)^\perp 
= \ker \,  \mP^\eps_\flat \, . $

\medskip

\noindent{Expand the function $ u^* \in L^{2*}_T $ in Fourier series 
and decompose the action of $ \mP^\eps_\flat $ in view of the Fourier
modes}

\medskip

$ u^*(t,x,\theta) = \sum_{k \in \ZZ_*} \, \u_k (t,x) \ e^{i \, 
k \, \theta} \, , \qquad \mP^\eps_\flat \, u^* = \sum_{k \in \ZZ_*} \, 
\mP^\eps_{\flat k} \, \u_k (t,x) \ e^{i \, k \, \theta}\, . $ 

\medskip

\noindent{Simple computations indicate that}

\medskip

$ \mP^\eps_{\flat k} \,\u_k \, := \, e^{- \, i \, \eps^{-1} \, k \, 
\varphi^\eps_\flat} \ \Pi (D_x) \ \bigl( e^{i \, \eps^{-1} \, k \, 
\varphi^\eps_\flat} \, \u_k \bigr) \, . $

\medskip

\noindent{The following result explains why the projector 
$ \mP^\eps_\flat $ is replaced by $ \Pi_0 $ when performing 
the BKW calculus.}

\begin{lem} \label{projeler} $ \, $

\smallskip

\noindent{i) The family $ \{ \mP^\eps_\flat \}_\eps $ is in 
$ \mU \mL^0 $. We have $ [ \part_\theta ; \mP^\eps_\flat ] 
= 0 $ and}

\medskip

$ [ \md_{j,\eps} ; \mP^\eps_\flat ] = 0 \, , \qquad \forall \, 
j \in \{ 0 , \cdots , d \} \, . $

\medskip

\noindent{ii) The projector $ \Pi^\eps_\flat (t,x) $ is an 
approximation of $ \mP^\eps_\flat $ in the sense that}

\medskip

$ \bigl \{ \, \mP^\eps_\flat - \Pi^\eps_\flat \, \bigr \}_\eps 
\, \in \, \eps \ \mU \mL^{2+\frac{d}{2}} \, , \qquad \bigl \{ \, 
\mP^\eps_\flat \ (\id - \Pi^\eps_\flat) \, \bigr \}_\eps \, \in \, 
\eps \ \mU \mL^1 \, . $

\end{lem}

\noindent{\em \underline{Proo}f \underline{o}f \underline{the} 
\underline{Lemma} 4.\underline{2}.} Since $ \mP^\eps_\flat $ 
is a projector, we are sure that

\medskip

$ \parallel \mP^\eps_\flat \, u \parallel_{L^2_T} \ \leq \ 
\parallel u \parallel_{L^2_T} \, , \qquad \forall \, (\eps,u) 
\in \, ] 0, \eps_0 ] \times L^2_T \, . $

\medskip

\noindent{It shows that $ \{ \mP^\eps_\flat \}_\eps \in \mU \mL^0 $.}
Compute

\medskip

$ [ \md_{j,\eps} ; \mP^\eps_\flat ] \, u^* (t,x,\theta) \, 
= \, \sum_{k \in \ZZ_*} \, [ \, \eps \ \part_j + i \ k \ 
\part_j \varphi^\eps_\flat \, ; \, \mP^\eps_{\flat k} \, ] \, 
\u_k (t,x) \ e^{i \, k \, \theta} \, . $ 

\medskip

\noindent{Observe that}

\medskip

$ ( \eps \ \part_j + i \ k \ \part_j \varphi^\eps_\flat ) \
\mP^\eps_{\flat k} \, \u_k \, = \, e^{- \, i \, \eps^{-1} \, 
k \, \varphi^\eps_\flat} \ \Pi (D_x) \ \eps \, \part_j \ 
\bigl( e^{i \, \eps^{-1} \, k \, \varphi^\eps_\flat} \ 
\u_k ) $

\smallskip

\hfill $ = \, \mP^\eps_{\flat k} \ ( \eps \ \part_j + 
i \ k \ \part_j \varphi^\eps_\flat ) \, \u_k \, . \qquad 
\qquad \qquad \qquad \, $

\medskip

\noindent{All these informations give access to the first 
assertion i).} Now consider ii). The asymptotic expansion
formula for pseudodifferential operators say that for 
all $ \u_k $ in $ C^\infty_0 ( \RR^d_T ) $ we have
$$ \quad \forall \, (t,x) \in \RR^d_T \, , \quad \ 
\lim_{\eps \, \longrightarrow \, 0} \ \ \bigl \{ \,
( \mP^\eps_{\flat k} \, \u_k) (t,x) \, - \, \Pi 
\bigl( \nabla \varphi^\eps_\flat (t,x) \bigr) \, 
\u_k (t,x) \, \bigr \} \, = \, 0 \, . $$
Since $ \Pi^\eps_\flat = \Pi ( \nabla \varphi^\eps_\flat ) $,
it indicates that $ \mP^\eps_\flat $ is close to 
$ \Pi^\eps_\flat $. We have to make this information 
more precise. To this end, proceed to the decomposition

\medskip

$ u^* \, = \, v^* + \eps \ \nabla p^* + \part_\theta p^* 
\times X^\eps_\flat \, , \qquad v^* = \mP^\eps_\flat \, 
u^* \, . $

\medskip

\noindent{We seek a solution $ ( v^*, p^*) $ of these
constraints such that}

\medskip

$ v^* \, = \, \Pi^\eps_\flat \, u^* + \eps \ \tilde v^* \, , \qquad 
p^* \, = \, \parallel X^\eps_\flat \parallel^{-2} \ X^\eps_\flat \cdot
\part_\theta^{-1} u^* + \eps \ \tilde p^* \, . $

\medskip

\noindent{After substitution, we find the relation}

\medskip

$ - \, \nabla \bigl( \parallel X^\eps_\flat \parallel^{-2} \
X^\eps_\flat \cdot \part_\theta^{-1} u^* \bigr) \, = \, 
\tilde v^* + \eps \ \nabla \tilde p^* + \part_\theta \tilde 
p^* \times X^\eps_\flat $

\medskip

\noindent{which must be completed by the condition}

\medskip

$ - \, \Div \, ( \Pi^\eps_\flat \, u^* ) \, = \, \eps \ 
\Div \, \tilde v^* + X^\eps_\flat \cdot \part_\theta 
\tilde v^* \, . $

\medskip

\noindent{It follows that}

\medskip

$ \tilde v^* \, = \, - \, \mP^\eps_\flat \, \bigl \lbrack
\nabla ( \parallel X^\eps_\flat \parallel^{-2} \ X^\eps_\flat 
\cdot \part_\theta^{-1} u^* ) \bigr \rbrack + ( \mP^\eps_\flat 
- \id ) \ \mr \mi \Dive^\eps_\flat \ \bigl( \Div \, ( \Pi^\eps_\flat 
\, u^* ) \bigr) \, . $ 

\medskip

\noindent{In view of this relation, the point ii) becomes clear.
\hfill $ \Diamond $

\medskip

\noindent{Consider the Cauchy problem}

\medskip

$ \md_{0,\eps} u^* + \eps^{-1} \ \Grade^\eps_\flat \, p^* \, 
= \, f^* \, , \qquad \Dive^\eps_\flat \, u^* = 0 \, , \qquad 
u^* (0,\cdot) = h^* (\cdot) $

\medskip

\noindent{with data $ f^* \in L^{2*}_T $ and $ h^* \in L^{2*} $.} 
Compose on the left with $ \mP^\eps_\flat $. It yields 

\medskip

$ \md_{0,\eps} u^* = \mP^\eps_\flat \, f^* + [\md_{0,\eps} ; 
\mP^\eps_\flat ] \, u^* \, , \qquad u^* (0,\cdot) = \mP^\eps_\flat 
\, h^* (\cdot) \, . $

\medskip

\noindent{The Cauchy problem can be solved in two steps.} 
First extract $ u^* $ from the above equation. Then recover 
$ p^* $ from the remaining relations.

\bigskip

\noindent{$ \bullet $ {\bf Proof of the Proposition 
\ref{complement}.}} It remains to absorb the term $ \tilde 
g^\eps_\flat $. Use the decomposition

\medskip

$ g^\eps_\flat = \langle \tilde g^\eps_\flat \rangle + 
\tilde g^{\eps *}_\flat \, , \qquad \langle \tilde g^\eps_\flat 
\rangle \in \text{Im} \, (\Div) \, , \qquad \tilde g^{\eps *}_\flat 
\in \text{Im} \, (\Dive^\eps_\flat) \, . $

\medskip

\noindent{It suffices to choose}

\medskip

$ c \u^\eps_\flat \, := \, - \, \text{ridiv} \, \langle \tilde 
g^\eps_\flat \rangle - \mr \mi \Dive^\eps_\flat \, \tilde 
g^{\eps*}_\flat = \bigcirc (\eps^{\frac{N+1}{l}}) \, . $

\section{Stability of strong oscillations} $ \, $

\vskip -5mm

\noindent{The case of turbulent regimes ($ l \geq 3 $) 
will not be undertaken here.} From now on, fix $ l = 2 $ 
and $ N \gg ( 6 + d ) $.} Consider the Cauchy problem
\begin{equation} \label{CPdep}
\left \{ \begin{array}{l}
\part_t \u^\eps + ( \u^\eps \cdot \nabla ) \, \u^\eps + 
\nabla \p^\eps \, = \, \nu \ E^{\eps l}_{\flat 0} ( \part )
\, \u^\eps \, , \qquad \Div \ \u^\eps \, = \, 0 \, ,
\quad \\
\u^\eps (0,x) \, = \, \u^\eps_\flat (0,x) \, . 
\end{array} \right.
\end{equation}
Let $ T_\eps $ be the upper bound of the $ T \geq 0 $ such 
that (\ref{CPdep}) has a solution $ \u^\eps \in \cW^0_T $. 
Classical results \cite{Chem} for fluid equations imply that 
$ T_\eps > 0 $. Our aim in this chapter 5 is to investigate 
the singular limit `$ \eps $ {\it goes to zero}'. Such an 
analysis must at least contain the two following parts.

\medskip

\noindent{a) An {\it existence} result for a time $ T_0 $ which is 
independent on the small parameter $ \eps \in \, ]0,\eps_0] $.} 
It is required that

\medskip

$ \inf \ \{ \, T_\eps \, ; \ \eps \in \, ]0,1] \, \} \, \geq \, 
T_0 \, > \, 0 \, . $

\medskip

\noindent{When $ \nu > 0 $, or when $ \nu = 0 $ and $ d = 2 $,
we know \cite{Chem}-\cite{Li} that $ T_\eps = + \infty $ so 
that $ T_0 = + \infty $.} When $ \nu = 0 $ and $ d \geq 3 $, 
nothing guarantees that $ T_0 > 0 $. To our knowledge, this
is an open question.

\medskip

\noindent{b) A {\it convergence} result.} The exact solution 
$ \u^\eps $ is not sure to remain close on the whole interval
$ [0,T_0] $ to the approximate solution $ \u^\eps_\flat $ 
given by the Theorem \ref{appBKW}. Proving estimates on 
$ \u^\eps - \u^\eps_\flat $ is a delicate matter.

\subsection{Various types of instabilities}

$ \bullet $ {\bf Obvious instabilities.} The obvious instabilities
are the mechanisms of amplifications which can be detected by 
looking directly at the formal expansions $ \u^\eps_\flat $. 
They imply the non linear instability of Euler equations. 
Indeed, fix any $ T > 0 $, any $ \u_0 \in \cW^\infty_T 
(\RR^d) $ which is solution of $ (\cE) $, and any 
$ \delta > 0 $.} Work on the balls

\medskip

$ B_0 ( \u_0 ; \delta ] \, := \, \bigl \{ \, \u \in L^2 \, ;
\ \parallel \u (\cdot) - \u_0 (0,\cdot) \parallel_{L^2 (\RR^d)} 
\, \leq \, \delta \, \bigr \} \, . $ 

\medskip

$ B_T ( \u_0 ; \delta ] := \, \bigl \{ \, \u \in L^2_T \, ;
\ \parallel \u - \u_0 \parallel_{L^2 ([0,T] \times \RR^d)} \, 
\leq \, \delta \, \bigr \} \, . $ 

\begin{prop} \label{euinstap} For all constant $ C > 0 $, there 
are small data 

\medskip

$ (\h,\tilde \h) \in \bigl( \, B_0 ( \u_0 ; \delta ] 
\cap H^\infty \, \bigr)^2 \, , \qquad (\f, \tilde \f) 
\in  \bigl( \, B_T ( \u_0 ; \delta ] \cap \cW^\infty_T 
\, \bigr)^2 $ 

\medskip

\noindent{so that the Cauchy problems}

\medskip

$ \part_t \u + ( \u \cdot \nabla ) \, \u + \nabla \p \, =
\, \f \, , \qquad \Div \ \u \, = \, 0 \, , \qquad \u (0 , 
\cdot) = \h(\cdot) \, , $

\smallskip

$ \part_t \tilde \u + ( \tilde \u \cdot \nabla ) \, \tilde \u 
+ \nabla \tilde \p \, = \, \tilde \f \, , \qquad \Div \ \tilde 
\u \, = \, 0 \, , \qquad \tilde \u (0 , \cdot) = \tilde \h (\cdot) \, , $

\medskip

\noindent{have solutions $ (\u , \tilde \u) \in B_T (\u_0;\delta]^2 $ 
and there is $ t \in \, ]0,T] $ such that}
\begin{equation} \label{contrsta}
\begin{array}{rl}
\ \parallel (\u - \tilde \u)(t, \cdot) \parallel_{L^2 (\RR^d)} \ 
\geq \, C & \! \! \! \bigl( \, \parallel \h - \tilde \h \parallel_{L^2 
(\RR^d)} \\
\ & + \, \hbox{$\int_0^t$} \, \parallel (\f - \tilde \f) (s,\cdot) 
\parallel_{L^2 (\RR^d)} \ ds \, \bigr) \, . 
\end{array} 
\end{equation}

\end{prop}

\medskip

\noindent{Inequalities as (\ref{contrsta}) are well-known.} In 
general \cite{CGM}-\cite{FSV}-\cite{Gr}, the demonstration is
achieved in two steps. First detect equilibria where instability
arises in the discrete spectrum. Then establish that linearized
instability implies non linear instability. The procedure we adopt 
below is different.} We just look at approximate solutions like  
$ \u^\eps_\flat $. It follows a more simple proof of (\ref{contrsta}).

\medskip

\noindent{\em \underline{Proo}f \underline{o}f \underline{the} 
\underline{Pro}p\underline{osition} \underline{\ref{euinstap}}.}  
Take $ l=2 $ and $ N \geq (8+d)  $. Consider two deals of initial data

\medskip

$ \tilde U^1_k (0,x,\theta) \, , \qquad \varphi^1_k (0,x) \, , 
\qquad 1 \leq k \leq N \, , $

\medskip

$ \tilde U^2_k (0,x,\theta) \, , \qquad \varphi^2_k (0,x) \, , 
\qquad 1 \leq k \leq N \, . $

\medskip

\noindent{Fix these expressions in the following way}

\medskip

$ \tilde U^1_1 (0,\cdot) \equiv \tilde U^2_1 (0,\cdot) \, , 
\quad \ \varphi^1_1 (0,\cdot) \equiv \varphi^2_1 (0,\cdot)
\equiv 0 \, , \quad \ \varphi^1_2 (0,\cdot) \equiv \varphi^2_2 
(0,\cdot) \equiv 0 \, . $

\medskip

\noindent{It implies that}

\medskip

$ \tilde U^1_1 (t,\cdot) \equiv \tilde U^2_1 (t,\cdot) \, , 
\qquad \varphi^1_1 (t,\cdot) \equiv \varphi^2_1 (t,\cdot) 
\equiv 0 \, , \qquad \forall \, t \in [0,T] \, . $

\medskip

\noindent{Adjust $ \tilde U^1_2 (0,\cdot) $ and $ \tilde 
U^2_2 (0,\cdot) $ so that}

\medskip

$ \part_t \, ( \varphi^1_2 - \varphi^2_2 ) (0,\cdot) \, = \, - \,
\nabla \varphi_0 \cdot \langle \tilde U^1_2 - \tilde U^2_2 \rangle
(0,\cdot) \, \not \equiv \, 0 \, . $

\medskip

\noindent{Therefore, we are sure to find some $ t > 0 $ 
such that $ ( \varphi^1_2 - \varphi^2_2 ) (t,\cdot) 
\not \equiv 0 $.} It follows that
\begin{equation} \label{difama} 
\begin{array} {ll}
U^1_1 (t,x,\theta) \! \! \! & = \, \tilde U^1_1 \bigl(
t, x, \theta + \varphi^1_2( t, x) \bigr) \\
\ & \not \equiv \, U^2_1 (t, x,\theta) \, = \, \tilde 
U^1_1 \bigl(t, x,\theta + \varphi^2_2 (t,x) \bigr) \, . 
\qquad \qquad \
\end{array}
\end{equation}
Note $ \u^{\eps 1}_\flat $ and $ \u^{\eps 2}_\flat $ 
the approximate solutions built with the profiles 
$ \{ U^1_k \}_k $ and $\{ U^2_k \}_k $. The associated 
error terms are $ \f^{\eps 1}_\flat $ and $ \f^{\eps 
2}_\flat $.

\medskip

\noindent{Proceed by contradiction.} Suppose that the 
Proposition \ref{euinstap} is wrong. Then, there is 
$ C > 0 $ and $ \eps_1 \in \, ]0,\eps_0 ] $ such 
that for all $ \eps \in \, ] 0 ,\eps_1 ] $, we have
$$ \ \begin{array}{rl}
\parallel (\u^{\eps 1}_\flat - \u^{\eps 2}_\flat)(t, \cdot) 
\parallel_{L^2 (\RR^d)} \ \leq \, C & \! \! \! \bigl( \, \parallel 
(\u^{\eps 1}_\flat - \u^{\eps 2}_\flat)(0, \cdot) \parallel_{
L^2 (\RR^d)} \\
\ & + \, \hbox{$ \int_0^t $} \, \parallel (\f^{\eps 1}_\flat - 
\f^{\eps 2}_\flat) (s,\cdot) \parallel_{L^2 (\RR^d)} \ d s \, 
\bigr) \, .  
\end{array} $$
Divide this inequality by $ \sqrt \eps $. By construction, 
we have

\medskip

$ \eps^{- \frac{1}{2}} \ \parallel (\u^{\eps 1}_\flat -
\u^{\eps 2}_\flat)(0, \cdot) \parallel_{L^2 (\RR^d)} \, = \, 
\bigcirc ( \sqrt \eps ) \, , $

\medskip

$ \eps^{- \frac{1}{2}} \ \parallel (\f^{\eps 1}_\flat -  
\f^{\eps 2}_\flat)(s, \cdot) \parallel_{L^2 (\RR^d)} \, 
= \, \bigcirc ( \sqrt \eps ) \, , \qquad \forall \, s
\in [0,t] \, . $

\medskip

$ \eps^{- \frac{1}{2}} \ \parallel (\u^{\eps 1}_\flat -  
\u^{\eps 2}_\flat)(t, \cdot) \parallel_{L^2 (\RR^d)} $

\smallskip

\hfill $ = \, \parallel (U^1_1 - U^2_1) (t,\cdot, \eps^{-1} \, 
\varphi^\eps_g (t, \cdot) \parallel_{L^2 (\RR^d)} \, + \, 
\bigcirc ( \sqrt \eps ) \, . \qquad $

\medskip

\noindent{It follows that}
$$ \ \lim_{ \eps \, \longrightarrow \, 0 } \ \  \eps^{- 
\frac{1}{2}} \ \parallel (\u^{\eps 1}_\flat - \u^{\eps 
2}_\flat)(t, \cdot) \parallel_{L^2 (\RR^d)} \, = \, 
\parallel (U^1_1 - U^2_1)(t, \cdot) \parallel_{L^2 
(\RR^d \times \TT)} \, = \, 0 $$
which is inconsistent with (\ref{difama}). \hfill $ \Box $

\bigskip

\noindent{{\it Remark 5.1.1:}} In the demonstration presented 
above, the amplification is due to $ \varphi_2 $ which is the 
principal term in the adjusting phase. The presence of $ \varphi_2 $
becomes efficient in comparison with the other effects when 

\medskip

$ \vert \, \tilde U^1_1 \bigl(t,x,\theta + \varphi^1_2 (t,x) \bigr) -  
\tilde U^2_1 \bigl(t,x,\theta + \varphi^2_2 (t,x) \bigr) \, \vert \, 
\sim \, c \ t \, \gg \, \sqrt \eps \, . $

\medskip

\noindent{This requires to wait a lapse of time bigger than 
$ \sqrt \eps $.} This delay can be reduced by adapting the 
above procedure to the cases $ l > 2 $. \hfill $ \triangle $

\medskip

\noindent{Obvious instabilities have an important consequence.}
To describe the related amplifications, it is necessary to introduce 
new quantities which correspond to the phase shifts. In other words, 
the only way to get $ L^2 \, - $estimates is to {\it blow up} the 
state variables. This principle is detailed in \cite{Che} in the 
case of compressible Euler equations.

\bigskip

\noindent{$ \bullet $ {\bf Hidden instabilities.}} Hidden 
instabilities are the amplifications which are not detected by 
the monophase description of the section 4. On the other hand, 
they can be revealed by a multiphase analysis. Introduce a second 
phase $ \psi_0 (t,x) \in \cW^\infty_T $ such that

\medskip

$ \part_t \psi_0 + (\u_0 \cdot \nabla) \, \psi_0 = 0 \, , 
\qquad \nabla \psi_0 \wedge \nabla \varphi_0 \not \equiv 0 $

\medskip

\noindent{and disturb the Cauchy data of (\ref{CPdep}) 
according to}

\medskip

$ \u^\eps (0,x) \, = \, \u^\eps_\flat (0,x) \, + \, \eps^{\frac
{M}{l}} \ U \bigl(x, \eps^{-1} \ \psi_0 (0,x) \bigr) \, , \qquad 
M \gg N \, . $

\medskip

\noindent{The small oscillations contained in the perturbation 
of size $ \eps^{\frac{M}{l}} $ are not always kept under control. 
They interact with $ \u^\eps_\flat $ and with themselves. They can 
be organized in such a way to affect the leading oscillation $ \u^\eps_\flat $. 
Concretely (see \cite{CGM}), we can adjust $ U $ and $ \psi_0 $ so that
there is a constant $ C > 0 $ and times $ t_\eps \in \, 
]0,T_\eps[ $ going to zero with $ \eps $ such that

\medskip

$ \parallel (\u^\eps -   \u^\eps_\flat) (t,\cdot) \parallel_{
L^2 (\RR^d)} \, \geq \, C \ \eps^{\frac{1}{2}} \, , \qquad 
\forall \, \eps \in \, ]0,\eps_0] \, . $

\medskip

\noindent{The power $ \eps^{\frac{M}{l}} $ at the time $ t=0 $ 
is turned into $ \eps^{\frac{1}{2}} $ at the time $ t=t_\eps $.}
Such amplifications occur whatever the selection of $ l \geq 2 $.
They imply minorations like (\ref{contrsta}). However, the underlying 
mechanisms are distinct from the preceding ones. They are implemented 
by oscillations which are transversal to $ \varphi_0 $ and whose 
wavelengths are $ \bigcirc(\eps ) $. They are cancelled by the 
addition of the anisotropic viscosity $ \nu \ E^{\eps l}_{\flat 0} $.

\subsection{Exact solutions}

$ \bullet $ {\bf Statement of the result.} The first information 
brought by the BKW construction is that mean values $ \bar U_k $ and 
oscillations $ U^*_k $ of the profiles $ U_k $ do not play the same
part. This fact is well illustrated by the rules of transformation 
(\ref{rulesoftr}). It means that we have to distinguish these quantities 
if we want to go further in the analysis. This can be done by involving 
the variables $ (t,x,\theta) $ that is by working at the level of 
(\ref{eclacpa}). To deal with $ (u^\eps , p^\eps) (t,x,\theta) $ 
instead of $ (\u^\eps , \p^\eps) (t,x) $ is usual in non linear 
geometric optics \cite{Sc}. It allows to mark the terms apt to 
induce instabilities.

\medskip

\noindent{Select some approximate solution $ (u^\eps_{(2,N)} , 
p^\eps_{(2,N)} ) $ with source term $ f^\eps_{(2,N)} $ given by 
the Proposition \ref{complement} and look at
\begin{equation} \label{eclacpasui}
\ \left \lbrace \begin{array} {l}
\md_{0,\eps} \, u^\eps + ( u^\eps \cdot \Grade^\eps_\flat) \, 
u^\eps + \Grade^\eps_\flat \, p^\eps \, = \, \nu \ \eps \ \Delta \, 
u^\eps \, , \qquad \Dive^\eps_\flat \, u^\eps = 0 \, , \\
u^\eps (0,x,\theta) \, = \, u^\eps_{(2,N)} (0,x,\theta) \, .
\end{array} \right.
\end{equation}

\begin{theo} \label{ciprin} Fix any integer $ N > d + 8 $. There 
is $ \eps_N \in \, ]0,1] $ and $ \nu_N > 0 $ such that for all 
$ \eps \in \, ]0,\eps_N] $ and for all $ \nu > \nu_N $ the Cauchy 
problem (\ref{eclacpasui}) has a unique solution $ (u^\eps , p^\eps) $ 
defined on the strip $ [0,T] \times \RR^d \times \TT $. Moreover

\medskip

$ \bigl \lbrace u^\eps -  u^\eps_{(2,N)}  \bigr \rbrace_\eps \, 
= \, \bigcirc (\eps^{\frac{N}{2}-d-4}) \, . $

\end{theo}

\smallskip

\noindent{\em \underline{Proo}f \underline{o}f \underline{the} 
\underline{Theorem} \underline{\ref{ciprin}}.} The system
(\ref{eclacpasui}) amounts to the same thing as
\begin{equation} \label{singudepmoyos}
\left \{ \begin{array}{l}
\part_t \bar u^\eps + (\bar u^\eps \cdot \nabla) \bar u^\eps 
+ \Div \, \langle u^{\eps*} \otimes u^{\eps*} \rangle + 
\nabla \bar p^\eps = \nu \ \Delta_x \, 
\bar u^\eps \, , \\ 
\part_{0,\eps} u^{\eps*} + (\bar u^\eps \cdot \Grade^\eps_\flat) \, 
u^{\eps*} + \eps \ (u^{\eps*}  \cdot \nabla ) \, 
\bar u^\eps \\
\qquad \quad \, + \, \left \lbrack ( u^{\eps *} \cdot 
\Grade^\eps_\flat) \, u^{\eps*} \right \rbrack^* + 
\Grade^\eps_\flat \, p^{\eps *} = \nu \ \eps \ 
\Delta \, u^{\eps *} \, , \qquad \ \\
\Div \, \bar u^\eps = \Dive^\eps_\flat \, u^{\eps *} = 0 \, . 
\end{array} \right.
\end{equation}
The equation (\ref{singudepmoyos}) is also equivalent to solve
the Cauchy problem
\begin{equation} \label{singudepmoyosmod}
\left \{ \begin{array}{l}
 P \, \part_t \bar u^\eps + P \, \left \lbrack (\bar u^\eps \cdot \nabla)
\bar u^\eps \right \rbrack + P \, \left \lbrack \Div \, \langle 
u^{\eps*} \otimes u^{\eps*} \rangle \right \rbrack = 
\nu \ \Delta_x \, \bar u^\eps \, , \qquad \\ 
\mP^\eps_\flat \, \part_{0,\eps} u^{\eps*} + \mP^\eps_\flat \, \left \lbrack 
(\bar u^\eps \cdot \Grade^\eps_\flat) \, u^{\eps*} \right 
\rbrack + \eps \ \mP^\eps_\flat \, \left \lbrack (u^{\eps*} 
\cdot \nabla ) \, \bar u^\eps \right \rbrack \\
\qquad \quad \, + \, \mP^\eps_\flat \, \left \lbrack 
( u^{\eps *} \cdot \Grade^\eps_\flat) \, u^{\eps*} 
\right \rbrack^* = \nu \ \eps \ \mP^\eps_\flat \ 
\Delta \, u^{\eps *} \, , 
\end{array} \right.
\end{equation}
associated with the compatible initial data

\medskip

$ \bar u^\eps (0, \cdot) = P \, \bar u^\eps_\flat (0, \cdot) \, ,
\qquad u^{\eps*} (0, \cdot) = \mP^\eps_\flat \, u^{\eps*}_\flat 
(0, \cdot) \, . $

\medskip

\noindent{$ \bullet $ {\bf Blow up.}} Introduce the new unknown 

\medskip

$ d^\eps \, = \, {}^t ( \bar d^\eps , d^{\eps*}) \, = \, {}^t 
( \, P \, \bar d^\eps \, , \, \mP^\eps_\flat \, d^{\eps*} \, ) $

\smallskip

$ \quad \ := \, 
\eps^{- \iota} \ \bigl( \, \eps^{-\frac{1}{l}} \ ( \bar u^\eps 
- \bar u^\eps_\flat) \, , \, ( u^{\eps *} - u^{\eps *}_\flat )
\bigr) \, , \qquad \flat = (2,N) \, . $ 

\medskip

\noindent{This transformation agrees with (\ref{rulesoftr}).} 
The weight $ \eps^{-\frac{1}{l}} $ in front of $ ( \bar u^\eps 
- \bar u^\eps_\flat)  $ induces a shift on the indice $ l $. 
Functions $ \bar U_l $ and $ U^*_{l-1} $ play now the same 
part related to the  amplifications. To write the equation 
on $ d^\eps $ in an abbreviated form, we need notations. 
Quasilinear terms 
$$ \ \left. \begin{array}{l}
\cL_{11}^\eps \, \bar d \, := \, P \, \bigl \lbrack \, 
(\bar u^\eps_\flat \cdot \nabla) \bar d \, \bigr \rbrack \, , 
\qquad \qquad \quad \ \\
\cL_{12}^\eps \, d^* \, := \, P \, \bigl \lbrack \, \Div 
\, \langle \eps^{- \frac{1}{2}} \ u^{\eps *}_\flat \otimes 
d^* \, + \, d^* \otimes \eps^{- \frac{1}{2}} \ u^{\eps 
*}_\flat \rangle \, \bigr \rbrack \, , \\
\cL_{21}^\eps \, \bar d \, := \, \eps^{\frac{1}{2}} \ 
\mP^\eps_\flat \ \bigl \lbrack \, ( u^{\eps *}_\flat 
\cdot \nabla) \bar d \, \bigr \rbrack \, , \\
\cL_{22}^\eps \, d^* \, := \, \mP^\eps_\flat \, \left 
\lbrack \, (\bar u^\eps_\flat \cdot \nabla ) d^* \, 
\right \rbrack \, + \, \eps^{-1} \ \mP^\eps_\flat \, \left 
\lbrack \, ( u^{\eps *}_\flat \cdot \Grade^\eps_\flat) \, 
d^* \, \right \rbrack^* \qquad \qquad \qquad \\
\qquad \qquad \ + \, \eps^{-1} \ \mP^\eps_\flat \, 
\left \lbrack \, (\partial_t \varphi^\eps_\flat + \bar 
u^\eps_\flat \cdot \nabla \varphi^\eps_\flat ) \ 
\partial_\theta d^* \, \right \rbrack \, . 
\end{array} \right. $$
Semilinear terms
$$  \ \left. \begin{array}{l}
A_{11}^\eps \, \bar d \, := \, P \, \bigl \lbrack \,
( \bar d \cdot \nabla) \bar u^\eps_\flat \, \bigr 
\rbrack \, , \qquad \\
A_{21}^\eps \, \bar d \, := \, \mP^\eps_\flat \, 
\bigl \lbrack \, ( \bar d \cdot \Grade^\eps_\flat) \, 
( \eps^{- \frac{1}{2}} \ u^{\eps *}_\flat ) \, \bigr 
\rbrack \, , \\
A_{22}^\eps \, d^* \, := \, \mP^\eps_\flat \, \left 
\lbrack \, ( d^* \cdot \nabla) \, \bar u^\eps_\flat \, 
\right \rbrack \, + \, \eps^{-1} \ \mP^\eps_\flat \, 
\left \lbrack \, ( d^* \cdot \Grade^\eps_\flat) \, 
u^{\eps*}_\flat \, \right \rbrack^* \, . \qquad 
\qquad \quad
\end{array} \right. $$
Small quadratic terms 
$$  \ \left. \begin{array}{l}
Q_1^\eps \, := \, \eps^{\frac{3}{2}} \ P \, \bigl 
\lbrack \, \Div \ ( \bar d \otimes \bar d ) \bigr 
\rbrack \, + \, \eps^{\frac{1}{2}} \ P \, \bigl 
\lbrack \, \Div \ \langle d^* \otimes d^* \rangle \, 
\bigr \rbrack \, , \\
Q_2^\eps \, := \, \eps^{\frac{1}{2}} \ \mP^\eps_\flat \, 
\left \lbrack \, (\bar d \cdot \Grade^\eps_\flat) \, d^* \,
\right \rbrack \, + \, \eps^{\frac{3}{2}} \ \mP^\eps_\flat \, 
\left \lbrack \, (d^* \cdot \nabla ) \, \bar d \, 
\right \rbrack \qquad \qquad \qquad \quad \ \ \\
\qquad \quad + \, \mP^\eps_\flat \, \left 
\lbrack \, ( d^* \cdot \Grade^\eps_\flat) \, d^* \, 
\right \rbrack^* \, .
\end{array} \right. $$
And error terms 

\medskip

$ er^\eps_1 \, := \, \eps^{-\iota-\frac{3}{2}} \ P \ 
\bar f^\eps_\flat \, , \qquad er^\eps_2 \, := \, 
\eps^{-\iota-1} \ \mP^\eps_\flat \ f^{\eps *}_\flat \, . $

\medskip
 
\noindent{With these conventions, the expression $ d^\eps $ 
is subjected to}
\begin{equation} \label{interface}
\ \left \{ \begin{array}{l}
P \, \partial_t \bar d^\eps + \cL^\eps_{11} \, \bar d^\eps 
+ \cL^\eps_{12} \, d^{\eps *} + A^\eps_{11} \, \bar d^\eps \\ 
\qquad \quad + \, \eps^{\iota - 1} \ Q^\eps_1 + er^\eps_1
\, = \, \nu \ P \ \Delta_x \bar d^\eps \, , \\
\mP^\eps_\flat \, \partial_t d^{\eps *} + \cL^\eps_{21} \, 
\bar d^\eps + \cL^\eps_{22} \, d^{\eps *} + A^\eps_{21} \, 
\bar d^\eps +  A^\eps_{22} \, d^{\eps *} \qquad \qquad 
\qquad \quad \ \\ 
\qquad \quad + \, \eps^{\iota - 1} \ Q^\eps_2 + er^\eps_2
\, = \, \nu \ \mP^\eps_\flat \ \Delta \, d^{\eps *} \, .
\end{array} \right.
\end{equation}

\bigskip

\noindent{Energy estimates are obtained at the level 
of (\ref{interface}).} Below, we just sketch the related 
arguments which are classical. 

\bigskip

\noindent{$ \bullet $ {\bf $ L^2 - \, $estimates for the 
linear problem.}} The linearized equations of Euler 
equations along the approximate solution $ u^\eps_\flat $ 
are obtained by removing $ Q^\eps_1 $ and $ Q^\eps_2 $ 
from (\ref{interface}). It yields a system which, at 
first sight, involves coefficients which are singular 
in $ \eps $. In fact, this is not the case. Let us 
explain why.

\medskip

\noindent{This is clear for $ \cL^\eps_{11} $, $ \cL^\eps_{21} $
and $ A^\eps_{11} $.}

\medskip

\noindent{Since $ u^{\eps*}_\flat = \bigcirc( \eps^{\frac{1}{l}} ) $, 
this is also true for $ \cL^\eps_{12} $ and $ A^\eps_{21} $.}

\medskip

\noindent{The contributions which in $ \cL^\eps_{22} $
have $ \eps^{-1} $ in factor give no trouble since}

\medskip

$ \part_t \varphi^\eps_\flat + \bar u^\eps_\flat \cdot \nabla 
\varphi^\eps_\flat = \bigcirc ( \eps^{\frac{N}{2}} ) = \bigcirc (
\eps^{d+4} ) \, , \qquad u^{\eps *}_\flat \cdot \nabla \varphi^\eps_\flat 
= v^{\eps *}_\flat = \bigcirc ( \eps^{1 + \frac{1}{l}} ) \, . $

\medskip

\noindent{Now, look at $ A^\eps_{22} $.} Recall that $ d^{\eps *} =
\mP^\eps_\flat \, d^{\eps *} $ which means that

\medskip

$ \eps^{-1} \ d^{\eps *} \cdot \nabla \varphi^\eps_\flat \, = \, 
- \, \Div \, d^{\eps *} \, . $

\medskip

\noindent{Therefore}

\medskip

$ \eps^{-1} \ \mP^\eps_\flat \, \left \lbrack \, ( d^{\eps *} 
\cdot \Grade^\eps_\flat) \, u^{\eps*}_\flat \, \right \rbrack^* \, 
= \, T^\eps (t,x,\nabla) \, d^{\eps *} \, , $

\medskip

\noindent{where $ T^\eps $ is some differential operator of 
order $ 1 $ with bounded coefficients.} 

\medskip

\noindent{Observe that these manipulations and the blow up procedure 
induce a {\it loss} of hyperbolicity. When $ \nu = 0 $, this is the 
source of hidden instabilities. When $ \nu \geq \nu_N > 0 $ with 
$ \nu_N $ large enough, this can be compensated by the viscosity. 
This is the key to $ L^2 - \, $estimates. 

\bigskip

\noindent{$ \bullet $ {\bf The non linear problem and higher 
order estimates.}} Let $ \sigma $ be the smaller integer such
that $ \sigma \geq \frac{d+3}{2} $. If the life span $ T_\eps $
of the exact solution $ u^\eps $ is finite, we must have
$$ \lim_{t \, \longrightarrow \, T_\eps} \ \ \parallel u^\eps 
(t,\cdot) \parallel_{H^\sigma} \, = \, + \infty \, . \qquad 
\qquad \qquad \qquad \qquad \qquad \qquad \qquad \qquad $$
Thus, the Theorem \ref{ciprin} is a consequence of the following 
majoration

\medskip 

$ \sup \ \bigl \{ \, \parallel u^\eps (t,\cdot) \parallel_{
H^\sigma} \, ; \ t \in [0 , \min \, (T_\eps,T) ] \, \bigr \} 
\, \leq \, C \, < \, \infty \, . $

\medskip

\noindent{Consider the set}

\medskip

$ \cZ_\eps \, := \, \bigl \{ \, \md_{0,\eps} \, , \, \cdots \, , 
\, \md_{d,\eps} \, , \, \partial_\theta \, \bigr \} \, . $

\medskip

\noindent{Extract the operators}

\medskip

$ \cZ^k_\eps \, := \, \cZ_1 \circ \, \cdots \, \circ \cZ_k \, , 
\qquad \cZ_j \in \cZ_\eps \, , \qquad k \leq \sigma \, . $

\medskip

\noindent{It suffices to show that}
$$ \ \max_{\, 0 \leq k \leq \sigma} \ \ \sup \ \bigl \{ \, 
\parallel \eps^{-k} \ \cZ^k_\eps \ u^\eps (t,\cdot) 
\parallel_{L^2} \, ; \ t \in [0 , \min \, (T_\eps,T) ] \, 
\bigr \} \, \leq \, C \, < \, \infty \, . \ $$
Pick some $ \cZ^k_\eps $ with $ k \leq \sigma $. Apply 
$ \cZ^k_\eps $ on the left of (\ref{interface}). Use 
the point i) of Lemma \ref{projeler} to pass through 
$ \mP^\eps_\flat $. Then, observe that the commutator 
of two vector fields in $ \cZ^\eps $ is a linear 
combination of elements of $ \cZ^\eps $ with 
coefficients in $ C^\infty $. Thus, we get an 
equation on $ \cZ^k_\eps \, d^{\eps *} $. 

\medskip

\noindent{The linear part is managed as in the preceding
paragraph.} Take $ \iota = 1 $. The contributions due to 
$ Q^\eps_1 $ and $ Q^\eps_2 $ are controled by way of the 
a priori estimate and the viscosity. The condition on $ N $ 
is to make sure that 

\medskip

$ \frac{N}{2} - \iota - \frac{3}{2} - \sigma \, \geq \, 0 \, . $

\medskip

\noindent{Thereby, the contributions brought by the error 
terms $ er^\eps_1 $ and $ er^\eps_2 $ remain bounded in 
the procedure.}

\bigskip

\bibliographystyle{alpha}

\end{document}